\documentclass[a4paper,12pt]{article}
\include{epsf}
\usepackage{latexsym}
\usepackage{amssymb}
\usepackage{amsmath}
\usepackage{amsxtra}
\usepackage{mathabx}
\usepackage[dvips]{graphicx}

\newlength{\cqfd}
\def\R{\mathbb{R}}
\def\C{\mathbb{C}}
\def\N{\mathbb{N}}
\def\Z{\mathbb{Z}}

\newtheorem{theorem}{Theorem}

\newtheorem{proposition}{Proposition}

\newtheorem{lemma}{Lemma}

\newcommand{\bq}{\overline{q}}
\newcommand{\ks}{k_\star}
\newcommand{\bqs}{\bq_\star}
\newcommand{\Ls}{L_\star}
\newcommand{\oms}{\omega_\star}
\newcommand{\cs}{c_\star}

\newcommand{\Hs}{H_\star}
\newcommand{\Qs}{Q_\star}
\newcommand{\Ms}{M_\star}

\newcommand{\Had}{H^{\txt{ad}}}
\newcommand{\Hads}{H^{\txt{ad}}_\star}

\newcommand{\GG}{\mathcal{G}}
\newcommand{\GGs}{\mathcal{G}_\star}
\newcommand{\LL}{\mathcal{L}}
\newcommand{\LLs}{\mathcal{L}_\star}
\newcommand{\KK}{\mathcal{K}}
\newcommand{\KKs}{\mathcal{K}_\star}
\newcommand{\LLad}{\mathcal{L}^{\txt{ad}}}

\newcommand{\ps}{p_\star}
\newcommand{\Pis}{\Pi_\star}

\newcommand{\txt}[1]{\textrm{#1}}
\renewcommand{\d}{\txt{d}}
\newcommand{\un}{\mathbf{1}}

\newcommand{\Ak}{A^k}
\newcommand{\AAk}{\widetilde{A}^k}
\newcommand{\Aq}{A^{\bq}}
\newcommand{\Aom}{A^{\omega}}
\newcommand{\Aks}{A_\star^k}
\newcommand{\AAks}{\widetilde{A}_\star^k}
\newcommand{\Aqs}{A_\star^{\bq}}
\newcommand{\Aoms}{A_\star^{\omega}}

\newcommand{\Bs}{B_\star}
\newcommand{\BTs}{B_{T,\star}}
\newcommand{\BXs}{B_{X,\star}}
\newcommand{\BBXs}{\widetilde{B}_{X,\star}}

\newcommand{\Xphi}{X^\varphi}
\newcommand{\Yphi}{Y^\varphi}
\newcommand{\hh}{\widetilde{h}}
\newcommand{\qq}{\widetilde{q}}

\newcommand{\BWT}{\mathcal{B}^T}
\newcommand{\BWX}{\mathcal{B}^X}
\newcommand{\BBWT}{\widetilde{\mathcal{B}}^T}
\newcommand{\BBWX}{\widetilde{\mathcal{B}}^X}

\newcommand{\As}{\mathcal{A}_\star}
\newcommand{\Asc}{\widecheck{\mathcal{A}}_\star}
\newcommand{\Ac}{\mathcal{A}^c}
\newcommand{\Acf}{\widehat{\mathcal{A}}^c}
\newcommand{\Qcc}{\widecheck{\mathcal{Q}}^c}
\newcommand{\PcMFc}{\widecheck{\mathcal{P}}^{c}_{\txt{mf}}}
\newcommand{\PcMF}{\mathcal{P}^{c}_{\txt{mf}}}
\newcommand{\PcFSc}{\widecheck{\mathcal{P}}^{c}_{\txt{fs}}}
\newcommand{\PcFS}{\mathcal{P}^{c}_{\txt{fs}}}
\newcommand{\Pcc}{\widecheck{\mathcal{P}}^{c}}

\newcommand{\PsMFc}{\widecheck{\mathcal{P}}^{s}_{\txt{mf}}}
\newcommand{\PsMF}{\mathcal{P}^{s}_{\txt{mf}}}
\newcommand{\PsFSc}{\widecheck{\mathcal{P}}^{s}_{\txt{fs}}}
\newcommand{\PsFS}{\mathcal{P}^{s}_{\txt{fs}}}
\newcommand{\Psc}{\widecheck{\mathcal{P}}^{s}}
\newcommand{\Ps}{\mathcal{P}^{s}}

\newcommand{\V}{\mathcal{V}}
\newcommand{\NL}{\mathcal{N}}
\newcommand{\NLc}{\mathcal{N}^{\txt{c}}}
\newcommand{\NLs}{\mathcal{N}^{\txt{s}}}
\newcommand{\NNLc}{\widetilde{\mathcal{N}}^{\txt{c}}}
\newcommand{\Lc}{\widecheck{\Lambda}}
\newcommand{\fc}[1]{\widecheck{#1}\sphat}
\newcommand{\cf}[1]{\widehat{#1}\spcheck}
\newcommand{\vc}{v^{\txt{c}}}
\newcommand{\vs}{v^{\txt{s}}}
\newcommand{\vceps}{v^{\txt{c}}_{\varepsilon}}
\newcommand{\vseps}{v^{\txt{s}}_{\varepsilon}}
\newcommand{\rceps}{r^{\txt{c}}_{\varepsilon}}
\newcommand{\rseps}{r^{\txt{s}}_{\varepsilon}}
\newcommand{\Vc}{V^{\txt{c}}}
\newcommand{\Vs}{V^{\txt{s}}}
\newcommand{\Vceps}{V^{\txt{c}}_{\varepsilon}}
\newcommand{\Vseps}{V^{\txt{s}}_{\varepsilon}}
\newcommand{\Resc}{\text{Res}_{\txt{c}}^\varepsilon}
\newcommand{\Ress}{\text{Res}_{\txt{s}}^\varepsilon}
\newcommand{\Nceps}{\mathcal{N}^{\txt{c}}_{\varepsilon}}
\newcommand{\Nseps}{\mathcal{N}^{\txt{s}}_{\varepsilon}}

\newcommand{\xnorm}[3]{\|#1\|_{\mathcal{X}^{#2}_{#3}}}
\newcommand{\ynorm}[3]{\|#1\|_{\mathcal{Y}^{#2}_{#3}}}
\newcommand{\hulnorm}[2]{\|#1\|_{H^{#2}_{\txt{ul}}}}

\title{\bf Whitham's equations for modulated roll-waves in shallow flows}

\begin{document}

\maketitle 
\begin{center}
{\large  Pascal Noble \footnote{ Universit\'e de Lyon, Universit\'e Lyon 1,
Institut Camille Jordan, UMR CNRS 5208, 43 bd du 11 novembre 1918,
F - 69622 Villeurbanne Cedex, France; noble@math.univ-lyon1.fr:
Research of P.N. was partially supported by French ANR project
no. ANR-09-JCJC-0103-01} L.Miguel Rodrigues  \footnote{Universit\'e de Lyon, Universit\'e Lyon 1,
Institut Camille Jordan, UMR CNRS 5208, 43 bd du 11 novembre 1918,
F - 69622 Villeurbanne Cedex, France; rodrigues@math.univ-lyon1.fr}}
\end{center}

\vspace{0.2cm}



$ $

\noindent
{\bf 2000 MR Subject Classification}: 35Q35, 35B10, 35A35, 35B35, 35P10, 35B27.


\vspace {0.5cm}
\begin{center}
{\bf Abstract.}
\end{center}

This paper is concerned with the detailed behaviour of roll-waves undergoing a low-frequency perturbation. We first derive the so-called \emph{Whitham's averaged modulation equations} and relate the well-posedness of this set of equations to the spectral stability problem in the small Floquet-number limit. We then fully validate such a system and in particular, we are able to construct solutions to the shallow water equations in the neighbourhood of modulated roll-waves profiles that exist for asymptotically large time.

\section{Introduction}

Our goal is to perform a two-scale analysis of waves in shallow flows, the fast scale being locally generated by (unmodulated) periodic travelling waves described by two parameters (for instance wavenumber and a discharge rate) and the low scale obeying an averaged system for these parameters. In studying such modulated waves our motivation is three-fold.

First, we are directly interested in a deeper understanding of such instabilities and, in this sense, this work is a piece of a wider program \cite{Ba_Jo_No_Ro_Zu_2, Ba_Jo_No_Ro_Zu_1,Ba_Jo_Ro_Zu,Jo_Zu_No,N1,N2}. Here we describe the motion of a shallow flow down an inclined ramp by the evolution of $(h,q)$, $h(t,x)\in\R^+$ being the fluid heigth and $q(t,x)\in\R$ the averaged horizontal momentum at time $t>0$ and in place $x\in\R$. The evolution we consider is governed by the Saint-Venant equations
\begin{equation}\label{sv}
\left\{\begin{array}{ll}
\displaystyle
\partial_t h+\partial_x q=0\\
\displaystyle
\partial_t q+\partial_x\Big(\frac{q^2}{h}+\frac{h^2}{2F^2}\Big)=h-\frac{q^2}{h^2}+\delta\partial_x^2q
\end{array}\right.
\end{equation}
taking into account viscosity ($\delta>0$, $\delta^{-1}$ being a Reynolds number), gravity ($F>0$ is a Froude number) and a turbulent (quadratic) friction along the bottom. The main (non obvious) physical flaw of this description lies probably in the form of the viscosity term $\delta\partial_x^2q$, which should be replaced with $\delta\partial_x(h\partial_x (q/h))$. But this restriction is purely motivated by writting convenience. The full paper would translate to the more physical case and actually even in this simplified case the only nonlinear system we solve is of quasilinear type. Roll-waves are then depicted as periodic travelling waves of system~\eqref{sv}, going down with a velocity larger than sound speed. Our first purpose is thus to investigate the behaviour of the solutions to these shallow water equations that are low-frequency perturbations of roll-waves.

The other two motivations consider roll-waves as representative of a wider class of periodic travelling waves. In the context of Lagrangian systems, Whitham explained how to derive an averaged system for the slow motion of the local parameters describing modulated periodic travelling waves \cite{Whitham} (see in particular Chapter $14$). In the following we will call such averaged systems, Whitham's systems. In the extension of the theory to a wider class of systems, and in a more mathematical way, Serre brought to light a direct relation between low-Floquet small eigenvalues of the original system linearized about a given periodic travelling wave and hyperbolicity of the Whitham's system linearized about corresponding paremeters \cite{Serre}. Since Whitham's systems are first order partial differential systems, this is precisely a relation between some spectral stability of the wave one wants to modulate and the well-posedness of the corresponding averaged system. Our goal is not only to extend this result to our situation (this is done in a straightforward way, see Lemma~\ref{evansserre}) but also to go one step further by associating not only eigenvalues but also eigenvectors (see Lemma~\ref{eigenvect1}). This is performed through a comparison of a spectral Fourier analysis of the averaged Whitham's system and a spectral Bloch analysis of the original Saint-Venant system. Note that this already provides us with at least a spectral validation of the Whitham's system.

Beyond the obvious interest of Lemma~\ref{eigenvect1} by itself, it is also a key-step for our third goal : to extend the work of Doelman, Sandstede, Scheel and Schneider \cite{D3S}, performed in the reaction-diffusion context, to system \eqref{sv}. Here we discuss only the inviscid part of \cite{D3S} and postpone to further work the viscous and shock (in parameters) parts. This extension will validate the Whitham's system at a nonlinear level, in the sense that to any solution to the Whitham's system, close to a given wave, one will associate a family of (higher-order) approximate solutions, with a modulated profile coincinding with the Whitham's solution at the linear level (see Proposition~\ref{approximate}), that describes at high order a family of solutions to system~\eqref{sv} for asymptotically large time (see Theorem~\ref{main}). Such a kind of validation of averaged equations has been performed in other contexts \cite{D3S,Dull_Schneider_NLS}. We believe the main difference between these cases and ours lies in the fact that we handle an averaged \emph{system}, our critical modes are not easily separated and we need a careful spectral preparation before being able to follow the strategy in \cite{D3S}. Lemma~\ref{eigenvect1} is precisely intended to fill this gap.

Together with the usual first-order Whitham's system (see system \eqref{whitham_xt}) we have discussed up to now, we also introduce an averaged second-order system (see system \eqref{whitham2_xt}), both in the derivation and in the spectral parts of our paper. In doing so we intend both to illustrate the strength of the spectral study through Bloch-Fourier comparison (see Lemma~\ref{eigenvect2}), and to prepare further work extending the viscous part of \cite{D3S} to our context. 

Our paper is organized as follows. In the first section, we set notations for the rest of this work. In the second one, we derive formally modulation systems and explain how to compute in the low-frequency regime modulated approximate solutions \emph{up to any order} with respect to $\varepsilon$ where $\varepsilon^{-1}$ is the characteristic wavelength of perturbations. Although these formal approximate solutions are not directly related to the ones we justify in Proposition~\ref{approximate}, their construction shed some light on the proof of Proposition~\ref{approximate}. Then we perform our spectral analysis. Afterwards, assuming the needed hyperbolicity of the Whitham's system, we will then provide a mathematical justification of this system in the spirit of what was done in \cite{D3S} in the reaction-diffusion framework. Finally we explain what are the main flaws of this justification and what may be expected from a detailed study of the second-order modulation system, postponed to further work.

$ $

\noindent{\bf Acknowledgement:} The authors warmly thank Arjen Doelman and Guido Schneider for kind enlightment about \cite{D3S}.

\section{Set-up}

\subsection{Existence of roll-waves}

We start recalling some properties of the set of periodic travelling-wave solutions to \eqref{sv}.

We search for a periodic travelling wave in the form 
$$
\displaystyle
(h,q)(x,t)=(H,Q)(\omega t+kx),
$$
\noindent
with $H,Q$ $1$-periodic functions and $\omega$, $k$ real numbers. Then $\omega$ is the time frequency, $k$ the wavenumber and $(H,Q)$ should satisfy the ordinary differential system
\begin{equation}\label{prof_sys}
\displaystyle
(\omega H+kQ)'=0,\quad \omega Q'+k\big(\frac{Q^2}{H}+\frac{H^2}{2F^2}\big)'=H-\frac{Q^2}{H^2}+\delta k^2 Q''.
\end{equation}
\noindent
Setting $\displaystyle c=-\omega/k$, the wave velocity, and integrating the first equation of \eqref{prof_sys} as $cH-Q=\bq$ yields the second order differential equation
\begin{equation}\label{profil}
\displaystyle
\delta k^2 c H''-k\left(\frac{\bq^2}{H}+\frac{H^2}{2F^2}\right)'+H-\left(c-\frac{\bq}{H}\right)^2=0.
\end{equation}

A first result, due to Dressler, yields the existence of {\it inviscid} roll-waves, which are necessarily discontinuous with Lax shocks as discontinuities.
\begin{proposition}[Dressler, $\delta=0$]\label{Dressler}
Let $F>2$ and $(k,\bq)\in\R^\star_+\times\R^\star_+$ fixed. Then there exists  a unique wavespeed $c=c^*(\bq)$, given by
$$
\displaystyle
c^*(\bq)=\bq^{\frac{1}{3}}\big(F^{\frac{1}{3}}+F^{-\frac{2}{3}}\big),
$$
such that there is a $1$-periodic  solution $H$ to (\ref{profil})$_{\delta=0}$; and this solution is unique (up to a translation). Moreover, the roll-wave can be alternatively parametrized by $(h_+,\bq)$ with $h_+$ the non dimensional minimum fluid height satisfying
$$
h_+>\frac{1}{F}+\frac{1}{2F^2}\big(1+\sqrt{1+4F}\big)=:h_m
$$
in such a way that 
$$
\begin{array}{ll}
\displaystyle
\frac{1}{k}=\frac{H_c}{F^2}\int_{h_+}^{h_-}\frac{h^2+h+1}{h^2-(F^{-2}+2F^{-1})h+F^{-2}}dh=:\frac{H_c}{F^2}\mathcal{P}(h_+),\\
\displaystyle
\int_0^1 H=\frac{H_c}{\mathcal{F}(h_+)}\int_{h_+}^{h_-}\frac{h^2+h+1}{h^2-(F^{-2}+2F^{-1})h+F^{-2}}h\,dh=:H_c\frac{\mathcal{Q}(h_+)}{\mathcal{P}(h_+)}.
\end{array}
$$
with $H_c=(\bq F)^{\frac{2}{3}}$ and $h_-$ uniquely determined by the Rankine Hugoniot jump condition: $\displaystyle h_-+h_+=\frac{2}{h_+h_-}$.
\end{proposition} 

When $\delta>0$ is fixed, one can prove the existence of \emph{small amplitude} continuous periodic travelling waves through a Hopf bifurcation argument. The existence result is even better as $\delta\to 0$ where the existence of large amplitude roll-waves, close to Dressler's roll-waves, is proved \cite{Haer, Noble}.

\begin{proposition}\label{prop2}
Let $F>2$ and $\bq>0$.
Then for any wavespeed $c$ such that $c^*(\bq)<c<c_{h}=c^*(\bq)+\mathcal{O}(\sqrt{\delta})$, there is a unique $k(c,\bq)$ such that there is a $1$-periodic solution $H$ to (\ref{profil})$_{\delta}$; and this solution is unique up to translation.\\
Moreover, for any fixed $\delta$ sufficiently small,  $\lim_{c\to c_h} k_\delta(c,\bq)=0$ (the continuous roll-wave converges to a solitary wave).\\
Alternatively, for any fixed $\delta$ sufficiently small, roll-waves can be para\-metrized by $(k,\bq)$; and $\lim_{\delta\to 0}c_{\delta}(k,\bq)=c^*(\bq)$ (the continuous roll-wave converges to a Dressler roll-wave as $\delta\to0$).
\end{proposition} 

In order to get a full idea of the bifurcation scenario, the reader is referred to \cite{nd,nd2,Hwang_Chang} where it is described. Yet note that there the viscosity term is really non physical and does not even provides us with the right jump condition in the small viscosity limit. However, we believe the scenario is the right one anyway and we check it numerically in \cite{Ba_Jo_No_Ro_Zu_2, Ba_Jo_No_Ro_Zu_1,Ba_Jo_Ro_Zu} (with a Lagrangian formulation).

We could work under the regime describe by Proposition~\ref{prop2} but we rather choose to take as an assumption that we will work in domains where the solutions to profile equation~\eqref{profil}, identified when coinciding up to translation, are uniquely and smoothly parametrized by $(k,\bq)$. Of course this assumption is stronger than the mere full-rank assumption in \cite{Serre}. We develop consequences of this assumption at the linear level in the next subsection.

\subsection{Abstract set-up}

In this subsection, we set some abstract notations for the rest of the paper.

For the sake of simplicity, let us first rewrite \eqref{sv} as 
\begin{equation}\label{sv_a}
\left\{\begin{array}{ll}
\displaystyle
\partial_t h+\partial_x q=0,\\
\displaystyle
\partial_t q+\partial_xG(h,q)=S(h,q)+\delta \partial_x^2q
\end{array}\right.
\end{equation}
where we have denoted
\begin{equation}
\displaystyle
G(h,q)=\frac{q^2}{h}+\frac{h^2}{2F^2},\quad S(h,q)=h-\frac{q^2}{h^2}.
\end{equation}
When looking for periodic travelling-wave solutions $(h,q)$ to \eqref{sv_a} with wavenumber $k\in\R_+^\star$, and frequency $\omega\in\R$, one must find $1$-periodic solutions $(H,Q)$ to
\begin{equation}\label{prof_sys_a}
\displaystyle
(\omega H+kQ)'=0,\quad \omega Q'+k\big(G(H,Q)\big)'=S(H,Q)+\delta k^2 Q''
\end{equation}
and set $(h,q)(t,x)=(H,Q)(\omega\,t+k\,x)$. For writting convenience, from now on we denote
\begin{equation}\label{GG}
\GG(H,Q;k)\ = \delta k^2 Q''-k\big(G(H,Q)\big)'+S(H,Q)
\end{equation}
and for later use
\begin{eqnarray}\label{GG_dh}
\partial_h \GG(H,Q;k)[f]\ =\ -k\big(\partial_h G(H,Q) f\big)'+\partial_h S(H,Q) f\\
\label{GG_dq}
\partial_q \GG(H,Q;k)[f]\ =\ \delta k^2 f''-k\big(\partial_q G(H,Q) f\big)'+\partial_q S(H,Q) f
\end{eqnarray}
and
\begin{equation}\label{GG_dk}
\partial_k \GG(H,Q;k)\ = \delta 2\,k Q''-\big(G(H,Q)\big)'\ .
\end{equation}

Denote also wave speed $c$ as $c=-\omega/k$. Then system \eqref{prof_sys_a} may be reduced to : there exists $\bq\in\R$ such that $c\,H-q=\bq$ and $H$ is a $1$-periodic solution to
\begin{equation}\label{profil_a}
\displaystyle
kc^2 H'\ +\ \GG(H,cH-\bq;k)\ =\ 0\ .
\end{equation}
We will work in a context where, once $(k,\bq)\in\R^\star\times\R$ fixed, profile equation \eqref{profil_a} possesses a $1$-periodic solution for one and only one speed $c(k,\bq)$ (therefore one frequency $\omega(k,\bq)=-kc(k,\bq)$) and, for $c=c(k,\bq)$, equation \eqref{profil_a} possesses a unique $1$-periodic solution, up to translation. Functions $c$ and $\omega$ are smooth and we also make a smooth choice of corresponding solution $H(k,\bq)=H(\,\cdot\,;k,\bq)$. Accordingly let us also denote $Q(k,\bq)=Q(\,\cdot\,;k,\bq)=c(k,\bq)H(\,\cdot\,;k,\bq)-\bq$. All these assumptions are of course generic and justified at least in the small-amplitude regime.

Now differentiating equation \eqref{profil_a} in $(k,\bq)$-variables leads to
\begin{equation}\label{LLOmega}
\begin{array}{rcl}
\LL(\ks,\bqs)\left[\d H(\ks,\bqs)\cdot(k,\bq)\right]&=&[\d \omega(\ks,\bqs)\cdot(k,\bq)]\ \Aom(\ks,\bqs)\\
&&+\ k\,\Ak(\ks,\bqs)\ +\ \bq\,\Aq(\ks,\bqs)
\end{array}
\end{equation}
or alternatively to
\begin{equation}\label{LLc}
\begin{array}{rcl}
\LL(\ks,\bqs)\left[\d H(\ks,\bqs)\cdot(k,\bq)\right]&=&-\ks[\d c(\ks,\bqs)\cdot(k,\bq)]\ \Aom(\ks,\bqs)\\
&& +\ k\,\AAk(\ks,\bqs)\ +\ \bq\,\Aq(\ks,\bqs)
\end{array}
\end{equation}
where $\LL(\ks,\bqs)$ is the linear operator associated to the linearization of equation \eqref{profil_a} around $H=H(\ks,\bqs)$, namely
\begin{equation}\label{LL}
\begin{array}{rcl}
\LL(\ks,\bqs)f&=&-\omega(\ks,\bqs)c(\ks,\bqs) f'\ +\ \partial_h \GG(H(\ks,\bqs),Q(\ks,\bqs);\ks)[f]\\
&& +\ c(\ks,\bqs)\partial_q \GG(H(\ks,\bqs),Q(\ks,\bqs);\ks)[f]\\
\end{array}
\end{equation}
and 
\begin{eqnarray}\label{Aom}
\Aom(\ks,\bqs)&=&-2\,\frac{\omega(\ks,\bqs)}{\ks}\,H(\ks,\bqs)'\\
\nonumber
&&+\ \frac{1}{\ks} \partial_q \GG(H(\ks,\bqs),Q(\ks,\bqs);\ks)[H(\ks,\bqs)]\\
\label{Aq}
\Aq(\ks,\bqs)&=&\partial_q \GG(H(\ks,\bqs),Q(\ks,\bqs);\ks)[\un]\\
\label{Ak}
\Ak(\ks,\bqs)&=&c(\ks,\bqs)^2\,H(\ks,\bqs)'\ -\ \partial_k \GG(H(\ks,\bqs),Q(\ks,\bqs);\ks)\\
\nonumber
&&-\ \frac{\omega(\ks,\bqs)}{\ks^2} \partial_q \GG(H(\ks,\bqs),Q(\ks,\bqs);\ks)[H(\ks,\bqs)]\\
\label{AAk}
\AAk(\ks,\bqs)&=&\Ak(\ks,\bqs)\ -\ c(\ks,\bqs)\Aom(\ks,\bqs)\\
\nonumber
&=&-c(\ks,\bqs)^2\,H(\ks,\bqs)'\ -\ \partial_k \GG(H(\ks,\bqs),Q(\ks,\bqs);\ks)\ .
\end{eqnarray}
Note that, associated to uniqueness of $c$ once $(\ks,\bqs)$ fixed comes the fact that $\Aom(\ks,\bqs)$ does not belong to the range of $\LL(\ks,\bqs)$ acting on $L^2_{\txt{per}}(\R)$, the space of $1$-periodic functions square-integrable on $]0,1[$.

Before going on with properties of $\LL(\ks,\bqs)$, let us choose to denote, for any function $f$,  $f_\star=f(\ks,\bqs)$. Likewise, in the modulational context, once fixed functions $(k_0,\bq_0)$ of variables $(X,T)$, for any function $f$, we will denote by $f_0$ corresponding values ; for instance
\begin{equation}\label{LL_0}
[\LL_0 f](y;X,T)\ =\ [\LL(k_0(X,T),\bq_0(X,T)) f(\,\cdot\,;X,T) ](y)\ .
\end{equation}

Coming back to $\LLs$, note that, due to translation invariance of equation \eqref{profil_a}, $\Hs'$ belongs to the kernel of $\LLs$. Moreover $\LLs$ is a Fredholm operator with index $0$. Let us denote $\LLad(\ks,\bqs)$ the formal adjoint operator of $\LLs$ and choose $\Had(\ks,\bqs)$ in its kernel such that
\begin{equation}\label{normal}
<\Had(\ks,\bqs);\Hs'>\ =\ 1\ 
\end{equation}
where $<\,\cdot\,;\,\cdot\,>$ is the scalar product on $L^2_{\txt{per}}(\R)$, namely, for $f,g\in L^2_{\txt{per}}(\R)$,
$$
<\,f\,;\,g\,>\ =\ \int_0^1 f\,g\ .
$$
Then, for any function $f$ belonging to $L^2_{\txt{per}}(\R)$, $f$ belongs to the range of $\LLs$ on $L^2_{\txt{per}}(\R)$ if and only if
\begin{equation}\label{cancel}
<\Hads;f>\ =\ 0\ .
\end{equation}
And we may define in a unique way an inverse 
\begin{equation*}
\LL^{-1}(\ks,\bqs)\ :\ \{\Hads\}^{\perp}\longrightarrow\{\Hads\}^{\perp}\ .
\end{equation*}
In the following, $\LLs^{-1}$ will always be meant to be defined in this way.

We will also need a projection on the range of $\LLs$. This may be done in a natural way by defining
\begin{equation}\label{projectors}
p(\ks,\bqs)f\ =\ \frac{<\Hads;f>}{<\Hads;\Aoms>}\quad\txt{and}\quad\Pi(\ks,\bqs)f\ =\ f-\ps(f)\Aoms\ .
\end{equation}
Then, for any $f$ in $L^2_{\txt{per}}(\R)$, equation 
\begin{equation*}
\LLs h\ =\ w \Aoms\ +\ f
\end{equation*}
has a solution $h\in L^2_{\txt{per}}(\R)$ if and only if $w=-\ps (f)$, and in this case a solution may be defined in a unique way by $h=\LLs^{-1}\Pis f$, and any solution may be written in a unique way
\begin{equation*}
h\ =\ \LLs^{-1}\Pis f\ +\ \alpha \Hs'
\end{equation*}
with $\alpha\in\R$.

More generally note that, for any $f$ in $L^2_{\txt{per}}(\R)$, equation
\begin{equation*}
\LLs h\ =\ k\Aks\ +\ \bq\Aqs\ +\ w \Aoms\ +\ f
\end{equation*}
has a solution $h\in L^2_{\txt{per}}(\R)$ if and only if
\begin{equation*}
w\ =\ \d\oms(k,s)\ -\ \ps(f)
\end{equation*}
and in this case a solution may be defined in a unique way by
\begin{equation*}
h\ =\ \d\Hs(k,\bq)\ +\ \LLs^{-1}\Pis f\ ,
\end{equation*}
and any solution may be written in a unique way
\begin{equation*}
h\ =\ \d\Hs(k,\bq)\ +\ \LLs^{-1}\Pis f\ +\ \alpha \Hs'
\end{equation*}
with $\alpha\in\R$.

Likewise, for any $f$ in $L^2_{\txt{per}}(\R)$, equation
\begin{equation*}
\LLs h\ =\ k\AAks\ +\ \bq\Aqs\ +\ w \Aoms\ +\ f
\end{equation*}
has a solution $h\in L^2_{\txt{per}}(\R)$ if and only if
\begin{equation*}
w\ =\ -\ks\d \cs(k,s)\ -\ \ps(f)
\end{equation*}
and in this case a solution may be defined in a unique way by
\begin{equation*}
h\ =\ \d\Hs(k,\bq)\ +\ \LLs^{-1}\Pis f\ ,
\end{equation*}
and any solution may be written in a unique way
\begin{equation*}
h\ =\ \d\Hs(k,\bq)\ +\ \LLs^{-1}\Pis f\ +\ \alpha \Hs'
\end{equation*}
with $\alpha\in\R$.

For later use, we will denote
\begin{equation}\label{KK}
\KK(\ks,\bqs)\ =\ \LLs^{-1}\Pis\ .
\end{equation}

As we will also have to solve equations like $h'\ =\ f$ for some $f$ in $L^2_{\txt{per}}(\R)$, it is worthwhile to note that this can be done in $L^2_{\txt{per}}(\R)$ if and only if 
\begin{equation*}
<f>\ :=\ \int_0^1 f\ =\ 0
\end{equation*}
and a solution may be defined in a unique way by
\begin{equation}\label{I}
h(y)\ =\ \int_0^y f\ =:\ I(f)(y)\ .
\end{equation}
Correspondingly, two more functions of $(k,\bq)$ will play a major role in the modulation analysis:
\begin{equation}\label{M,N}
M(k,\bq)\ =\ \int_0^1 H(y;k,\bq)\d y\ ,\quad N(k,\bq)\ =\ c(k,\bq)M(k,\bq)-\bq\ .
\end{equation}
Note that, whereas we made some choice in the parametrization of solutions $H(k,\bq)$, functions $M$ and $N$ do not depend on this choice.

Besides the assumption of smooth parametrization by $(k,\bq)$, we will almost always work with $(k,\bq)$ such that
\begin{equation}\label{evolution}
\partial_{\bq} M(k,\bq)\ \neq\ 0\ 
\end{equation}
which may be seen to hold both in small-amplitude and small-viscosity regimes. At last, let us say that we will also assume 
\begin{equation}\label{cparam}
\partial_k c(k,\bq)\ \neq\ 0\qquad\textrm{and\qquad}\partial_{\bq} c(k,\bq)\ \neq\ 0\ .
\end{equation}
The first part is motivated by the fact that sometimes we will switch to a $(c,\bq)$-parametrization and for this purpose we will need such an assumption but of course this assumption can be removed when we work with the usual $(k,\bq)$-parametrization. The second part is intended to simplify some discussions when solving linearized equations, and is not crucial.

\section{Formal derivation of modulation systems}

\subsection{First order Whitham's equations}

We follow the method proposed by Serre \cite{Serre} to derive Whitham's equations and first introduce rescaled variables $(X,T)=(\varepsilon x,\varepsilon t)$. This yields
\begin{equation}\label{sv_e}
\left\{\begin{array}{l}
\displaystyle
\partial_T h\ +\ \partial_X q\ =\ 0\ ,\\
\displaystyle
\partial_T q\ +\ \partial_XG(h,q)\ =\ \frac{S(h,q)}{\varepsilon}\ +\ \varepsilon\,\delta \partial_X^2q\ .
\end{array}\right.
\end{equation}
\noindent
We then search for an expansion of $(h,q)$ in the form
$$
\displaystyle
(h,q)=\sum_{i=0}^{\infty}\varepsilon^{i}(h^i,q^i)\left(\frac{\phi(X,T)}{\varepsilon}; X,T\right),
$$
with $(h^i,q^i)$ $1$-periodic in their first argument $y\in\mathbb{R}$. Identifying $\mathcal{O}(\frac{1}{\varepsilon})$ terms, one finds
$$\left\{
\begin{array}{ll}
\displaystyle
\partial_T \phi\,\partial_y h^0+\partial_X\phi\,\partial_y q^0=0,\\
\displaystyle
\partial_T \phi\,\partial_y q^0+\partial_X \phi\,\partial_y \big(G(h^0,q^0)\big)=S(h^0, q^0)+\delta(\partial_X\phi)^2\partial_y^2 q^0.
\end{array}\right.
$$
\noindent
Let us denote $k(X,T)=\partial_X \phi(X,T)$, $\omega(X,T)=\partial_T \phi(X,T)$ and $c=-\frac{\omega}{k}$ so that $k$ is the local wavenumber, $\omega$ the local frequency and $c$ the local wave speed. The previous system is solved by 
$$
\displaystyle
c(X,T)=c(k(X,T),\bq(X,T)),\quad \omega(X,T)=\omega(k(X,T),\bq(X,T))
$$ 
and  $(h^0,q^0)(y;X,T)=(H,Q)(y;k(X,T),\bq(X,T))$, where $\bq(X,T)$ is a local discharge rate. Note that $\phi$ may be recovered from $k$ and $\bq$ if and only if 
\begin{equation}\label{eq_mod_1}
\displaystyle
\partial_T k+\partial_X\left(k\,c(k,\bq)\right)=0.
\end{equation}

We further identify $\mathcal{O}(1)$ terms in \eqref{sv_e}. On the one hand, the mass conservation law yields
$$
\displaystyle
\partial_y (\omega(k,\bq)h^1-k q^1)=-\big(\partial_T h^0 +\partial_X q^0\big).
$$
This equation has a solution if and only if 
\begin{equation}\label{eq_mod_2}
\displaystyle
\partial_T <h^0>+\partial_X <q^0>\ =\ 0\ .
\end{equation}

Using the fact that $c\,h^0-q^0=\bq$, equation  \eqref{eq_mod_2} can be equivalently written as
\begin{equation}\label{eq_mod_3}
\displaystyle
\partial_T (M(k,\bq))+\partial_X(N(k,\bq))=0\ ,
\end{equation}
where $M$ and $N$ are defined in \eqref{M,N}. The system (\ref{eq_mod_1},\ref{eq_mod_3}) forms the system of {\it Whitham's equations}. Assumption \eqref{evolution} is equivalent to this system being of evolution-type for $(k,\bq)$.  This is a first order differential system of partial differential equations: in what follows , we will study the hyperbolicity of such a system and relate the hyperbolicity to the stability of roll-waves in the small wavenumber regime. We will use this set of equations to construct approximate solutions to the full shallow water system in the neighbourhood of roll-waves on asymptotically large time intervals.

\subsection{Higher order approximations}

In the following, we show how to construct a higher-order approximation of solutions to \eqref{sv_e}: for that purpose, we need also to expand the phase $\phi$ with respect to $\varepsilon$ just as in a classical WKB-type calculation. Solutions $(h,q)$ to \eqref{sv_e} are then expanded in the form
$$
\displaystyle
(h,q)=\sum_{i=0}^{\infty}\varepsilon^{i}(h^i,q^i)\left(\frac{\phi(X,T)}{\varepsilon},X,T\right),\quad \phi(X,T)=\sum_{j=0}^{\infty}\varepsilon^{j}\phi^{j}(X,T).
$$
with $(h^i,q^i)$ $1$-periodic in their first argument $y\in\mathbb{R}$. Identifying $\mathcal{O}(\varepsilon^{-1})$ still yields a differential system in the form
$$\left\{
\begin{array}{ll}
\displaystyle
\partial_T \phi^0\,\partial_y h^0+\partial_X \phi^0\,\partial_y q^0=0,\\
\displaystyle
\partial_T \phi^0\,\partial_y q^0+\partial_X \phi^0\,\partial_y \big(G(h^0,q^0)\big)=S(h^0, q^0)+\delta(\partial_X\phi^0)^2\partial_y^2 q^0.
\end{array}\right.
$$
As previously, we set
$$
\displaystyle
k_0=\partial_X \phi^0,\quad \omega^0=\partial_T \phi^0,\quad c^0=-\frac{\omega^0}{k_0}\ .
$$
There must be some $\bq_0$ such that $c^0h^0(y)-q^0(y)=\bq_0$. This yields
$$
\displaystyle
\omega^0\ =\ \omega(k_0,\bq_0),\quad c^0\ =\ c(k_0,\bq_0)
$$
and the profile equation in the $y$-variable
$$
\displaystyle
k_0 c_0^2 \partial_y h^0\ +\ \GG(h^0,c_0 h^0-\bq_0;k_0)\ =\ 0\ 
$$
solved by $h^0(y;X,T)=H(y;k_0,\bq_0)$. As already mentioned, the compatibility condition
$\displaystyle \partial_T\partial_X\phi^0=\partial_X\partial_T\phi^0$ yields an evolution equation for the local wavenumber $k_0$:
\begin{equation}
\displaystyle
\partial_T k_0+\partial_X\big(k_0c_0\big)=0\ .
\end{equation}
\noindent
Next, we identify $\mathcal{O}(1)$ terms. First, we consider the  mass conservation law:
$$
\displaystyle
\partial_y(c_0\,h^1-q^1)=\frac{1}{k_0}(\partial_T h^0+\partial_X q^0)+\frac{1}{k_0}(\partial_T \phi^1+c_0\partial_X\phi^1)\partial_y h^0.
$$
\noindent
Note that the expansion of $\phi$ with respect to $\varepsilon$ yields a new term that does not change the other compatibility condition
$$
\displaystyle
\partial_T M_0+\partial_X N_0\ =\ 0.
$$
\noindent
However the presence of this extra term will be necessary in order to compute an approximation to the next order of the solution to the Saint-Venant equations. Integrating this equation with respect to $y$ yields
\begin{equation}\label{eq_q1}
\displaystyle
c_0 h^1-q^1=\frac{1}{k_0}\big(\partial_T I(H_0)+\partial_X I(Q_0)\big)+(\omega_1+c_0k_1)\frac{H_0}{k_0}+\bq_1,
\end{equation}
\noindent
where $\bq_1$ is a constant of integration depending on $(X,T)$ which plays the role of a correction to the relative discharge rate $\bq_0$, $I$ was defined in \eqref{I} and
\begin{equation}
\omega_1\ =\ \partial_T \phi^1,\quad k_1\ =\ \partial_X \phi^1.
\end{equation} There remains to determine $h^1,\phi^1$. For that purpose, we consider the momentum equation to order $\mathcal{O}(1)$. This yields
\begin{equation}\label{eq_h3}
\LL_0 h^1\ =\ k_1\,\Ak_0\ +\ \bq_1\,\Aq_0\ +\ \omega_1\,\Aom_0\ +\ R_0
\end{equation}
with $R_0$ a function that depends only on $H_0, k_0,\bq_0$ and is defined as
{\setlength\arraycolsep{1pt}
\begin{eqnarray}\label{RO}
\displaystyle
R_0&=&\partial_T Q_0+\partial_X [G(H_0,Q_0)]-\delta\,\partial_X [k_0\partial_y Q_0] -\delta\,k_0\partial_y(\partial_X Q_0)\\
\displaystyle
&&c_0\partial_y[\partial_T I(H_0)+\partial_X I(Q_0)]\ +\ \frac{1}{k_0}\partial_q\GG(H_0,Q_0;k_0)[\partial_T I(H_0)+\partial_X I(Q_0)]\ .\nonumber
\end{eqnarray}}
\noindent

\noindent
Then equation \eqref{eq_h3} has a solution if and only if
\begin{equation}\label{eq_fi1}
\displaystyle
\omega_1\ -\ \d\omega_0 (k_1,\bq_1)\ =\ -\ p_0(R_0)\ 
\end{equation}
where $p$ is the function defined in \eqref{projectors}. Once again, in order to recover $\phi^1$ from $k_1$ and $\omega_1$, one must impose
\begin{equation}\label{eq_k1}
\displaystyle
\partial_T k_1\ -\ \partial_X[\d\omega_0 (k_1,\bq_1)]\ =\ -\partial_X[p_0(R_0)]\ .
\end{equation}

\noindent
With our choice for $\omega_1$, equation \eqref{eq_h3} may be solved by
\begin{equation}\label{def_h1}
h^1\ =\ \d H_0(k_1,\bq_1)\ +\ \KK_0 R_0,
\end{equation}
$\KK$ being the operator defined in \eqref{KK}. We still need an equation to couple with \eqref{eq_k1} in order to determine $(k_1,\bq_1)$: this is done by considering the mass conservation law to order $\mathcal{O}(\varepsilon)$. One finds
\begin{equation}\label{def_q2}
\displaystyle
\partial_y (c_0 h^2-q^2)=\frac{1}{k_0}\big((\partial_T h^1+\partial_X q^1)+\partial_y(\hdots)\big).
\end{equation}
\noindent
This equation has a solution provided that
\begin{equation}\label{cond_h1}
\displaystyle
\partial_T< h^1>+\partial_X< q^1>=0.
\end{equation}
Combining \eqref{eq_q1}, \eqref{eq_fi1} and \eqref{def_h1} yields
\begin{equation}\label{edp_dh1}
\displaystyle
\begin{array}{l}
\partial_T \left(\d M_0(k_1,\bq_1)\right)+\partial_X\left(\d N_0(k_1,\bq_1)\right)\\
\displaystyle
=\ -\partial_T<\KK_0 R_0>-\partial_X<c_0\KK_0 R_0>-\partial_X\left(\frac{M_0}{k_0}\,p_0(R_0)\right)\\
\displaystyle
\quad+\ \partial_X\left(\frac{1}{k_0}\left(\partial_T<I(H_0)>+\partial_X <I(Q_0)>\right)\right).
\end{array}
\end{equation}
\noindent

\noindent
The system (\ref{eq_k1}, \ref{edp_dh1}) on $(k_1,\bq_1)$ is linear hyperbolic provided that the nonlinear system (\ref{eq_mod_1}, \ref{eq_mod_3})  is hyperbolic and this former hyperbolicity is discussed in the next section. Under the assumption of hyperbolicity, we can compute $(k_1,\bq_1)$. Then, $h^1,q^1,\phi^1$ are fully determined by $(h^0,q^0,\phi^0)$. The next steps are done similarly: assume that we have determined $(h^i,q^i,\phi^i)_{i\leq n}$ and let us compute $(h^{n+1},q^{n+1},\phi^{n+1})$. First, we determine $q^{n+1}$ as a function of $h^{n+1},\phi^{n+1}$, $(h^i,q^i,\phi^i)_{i\leq n}$  and a $y$-independent $\bq_{n+1}$ with the help of an equation similar to \eqref{def_q2}. Indeed, this equation has a solution due to the compatibility condition similar to \eqref{cond_h1} with $(h^n,q^n)$ replacing $(h^1,q^1)$. Then, one inserts the expansion of $q^{n+1}$ into the momentum equation written at order $\mathcal{O}(\varepsilon^n)$. This yields an equation in the form
$$
\displaystyle
\LL_0 h^{n+1}\ =\ k_{n+1}\,\Ak_0\ +\ \bq_{n+1}\,\Aq_0\ +\ \omega_{n+1}\,\Aom_0+R_n, 
$$
\noindent
with $R_n$ a function of $(h^i,q^i,\phi^i)_{i\leq n}$ and $k_{n+1}=\partial_X\phi^{n+1}$, $\omega_{n+1}=\partial_T\phi^{n+1}$. This equation has a solution provided that
$$
\omega_{n+1}\ -\ \d\omega_0 (k_{n+1},\bq_{n+1})\ =\ -\ p_0(R_n)\ .
$$
This yields an evolution equation for $k_{n+1}$. Then one can solve the equation on $h^{n+1}$ by setting
$$
h^{n+1}\ =\ \d H_0(k_{n+1},\bq_{n+1})\ +\ \KK_0 R_n\ .
$$
One finishes the construction at the $(n+1)$-step by writing the compatibility condition
$$
\displaystyle
\partial_T< h^{n+1}>+\partial_X< q^{n+1}>=0\ 
$$ 
which provides $(k_{n+1},\bq_{n+1})$ with an equation similar to \eqref{edp_dh1} for $(k_1,\bq_1)$.
\noindent
As a consequence, we build $(k_{n+1},\bq_{n+1})$ by solving a linear hyperbolic system. This gives a full description of an approximate solution of order $\mathcal{O}(\varepsilon^{n+1})$ and completes the construction of approximate solution to \eqref{sv_e} up to {\it any} order with respect to $\varepsilon$.

In another section, we will rigorously validate this formal construction and use its principle to justify modulation system (\ref{eq_mod_1}, \ref{eq_mod_3}).

\subsection{Two modulation systems}

In this subsection we extract from the analysis of previous subsections two nonlinear systems for local wavenumber and local discharge rate $(k,\bq)$ that we think of great importance for the understanding of low-frequency perturbations of periodic travelling-wave solutions.

First, equations (\ref{eq_mod_1}, \ref{eq_mod_3}) form the Whitham's averaged system for the system \eqref{sv} :
\begin{equation}\label{whitham_XT}
\left\{
\begin{array}{rcl}
\partial_T k\ -\ \partial_X [\omega(k,\bq)]&=&0\\
\partial_T [M(k,\bq)]\ +\ \partial_X [N(k,\bq)]&=&0
\end{array}
\right.
\end{equation}
which stays unchanged when coming back to physical variables $(x,t)$ by setting $(k,\bq)(x,t)=(k,\bq)(\varepsilon X,\varepsilon T)$:
\begin{equation}\label{whitham_xt}
\left\{
\begin{array}{rcl}
\partial_t k\ -\ \partial_x [\omega(k,\bq)]&=&0\\
\partial_t [M(k,\bq)]\ +\ \partial_x[N(k,\bq)]&=&0\ .
\end{array}
\right.
\end{equation}
From \eqref{evolution} may be seen that this system is of evolution type. Note that, though this system does not contain second order terms, the presence of a viscosity in system \eqref{sv} was proeminent in order to obtain it. Moreover, whereas working with a realistic viscosity term would not have changed the form of the system (and the definition of $M$, $N$ and $\omega$), it may indeed have changed the actual values of $M$,$N$ and $\omega$.

Concerning the well-posedness of system \eqref{whitham_XT}, the best one can expect is that the system is hyperbolic and thus possesses local-in-time smooth solutions. In the next section, we will connect this hyperbolicity to the spectral stability of periodic travelling-wave solutions to system \eqref{sv} under low-frequency perturbations. Then, under this spectral stability assumption, we will validate system \eqref{whitham_XT} by proving that any smooth solution to \eqref{whitham_XT} describes at first order an $\varepsilon$-family of modulated \emph{true} solutions to \eqref{sv} on asymptotically large time.

To be able to go beyond the shock formation in finite time for solutions to \eqref{whitham_XT}, we derive now a second-order modulation system. This will be done by combining (\ref{eq_k1}, \ref{edp_dh1}) and (\ref{eq_mod_1}, \ref{eq_mod_3}). In order to do so, let us rewrite \eqref{RO} as 
\begin{eqnarray*}
R_0&=&B_{T,k}(k_0,\bq_0)\,\partial_Tk_0\ +\ B_{X,k}(k_0,\bq_0)\,\partial_Xk_0\\
&+&B_{T,\bq}(k_0,\bq_0)\,\partial_T\bq_0\ +\ B_{X,\bq}(k_0,\bq_0)\,\partial_X\bq_0
\end{eqnarray*}
where $B_T$ and $B_X$ are defined by
\begin{equation}\label{BT}
\begin{array}{rcl}
B_T(\ks,\bqs)(k,\bq)&=&\left[2\cs\d\Hs\,+\,\d\cs\,\Hs\right](k,\bq)\\
&+&\frac{1}{\ks}\partial_q\GG(\Hs,\Qs;\ks)[I(\d\Hs(k,\bq))]\ -\ q
\end{array}
\end{equation}
and
\begin{equation}\label{BX}
\begin{array}{rcl}
B_X(\ks,\bqs)(k,\bq)&=&\frac{1}{\ks}\partial_q\GG(\Hs,\Qs;\ks)[I(\Hs)]\d\cs(k,\bq)\\
&+&\left[\cs\,\Hs\,+\,\partial_q G(\Hs,\Qs)\Hs\,-\,2\delta\ks\Hs'\right]\d\cs(k,\bq)\\
&+&\left[\cs^2\d\Hs\,+\,\d G(\Hs,\Qs)(\d\Hs,\cs\d\Hs)\,-\,2\delta\ks\cs\d\Hs'\right](k,\bq)\\
&+&\frac{\cs}{\ks}\partial_q\GG(\Hs,\Qs;\ks)[I(\d\Hs(k,\bq))]\ -\ \delta\,\cs\Hs'\,k\\
&-&\left[\cs\,+\,\frac{1}{\ks}\partial_q\GG(\Hs,\Qs;\ks)[I(\un)]\,+\,\partial_q G(\Hs,\Qs)\right]\bq\ .
\end{array}
\end{equation}
Then, at order $\mathcal{O}(\varepsilon^2)$, the pair $(k,\bq)=(k_0,\bq_0)+\varepsilon(k_1,\bq_1)$ satisfies
\begin{equation}\label{whitham2_XT}
\left\{
\begin{array}{rcl}
\partial_T k\ -\ \partial_X [\omega(k,\bq)]&=&-\delta\,\partial_X\left(p(k,\bq)[B_T(k,\bq)\partial_T(k,\bq)]\right)\\
&-&\delta\,\partial_X\left(p(k,\bq)[B_X(k,\bq)\partial_X(k,\bq)]\right)\\
\partial_T [M(k,\bq)]\ +\ [\partial_X N(k,\bq)]&=&-\delta\,\partial_X\left(\frac{M(k,\bq)}{k}p(k,\bq)[B_T(k,\bq)\partial_T(k,\bq)]\right)\\
&-&\delta\,\partial_X\left(\frac{M(k,\bq)}{k}p(k,\bq)[B_X(k,\bq)\partial_X(k,\bq)]\right)\\
&+&\delta\,\partial_X\left(\frac{1}{k}<I(\d H(k,\bq)\partial_T(k,\bq))>\right)\\
&+&\delta\,\partial_X\left(\frac{1}{k}<I(\d Q(k,\bq)\partial_X(k,\bq))>\right)\\
&-&\delta\,\partial_T\left(<\KK(k,\bq)[B_T(k,\bq)\partial_T(k,\bq)]>\right)\\
&-&\delta\,\partial_T\left(<\KK(k,\bq)[B_X(k,\bq)\partial_X(k,\bq)]>\right)\\
&-&\delta\,\partial_X\left(c(k,\bq)<\KK(k,\bq)[B_T(k,\bq)\partial_T(k,\bq)]>\right)\\
&-&\delta\,\partial_X\left(c(k,\bq)<\KK(k,\bq)[B_X(k,\bq)\partial_X(k,\bq)]>\right)
\end{array}
\right.
\end{equation}
which in physical variables turns into
\begin{equation}\label{whitham2_xt}
\left\{
\begin{array}{rcl}
\partial_t k\ -\ \partial_x [\omega(k,\bq)]
&=&-\partial_x\left(p(k,\bq)[B_T(k,\bq)\partial_t(k,\bq)]\right)\\
&-&\partial_x\left(p(k,\bq)[B_X(k,\bq)\partial_x(k,\bq)]\right)\\
\partial_t [M(k,\bq)]\ +\ [\partial_x N(k,\bq)]
&=&-\partial_x\left(\frac{M(k,\bq)}{k}p(k,\bq)[B_T(k,\bq)\partial_t(k,\bq)]\right)\\
&-&\partial_x\left(\frac{M(k,\bq)}{k}p(k,\bq)[B_X(k,\bq)\partial_x(k,\bq)]\right)\\
&+&\partial_x\left(\frac{1}{k}<I(\d H(k,\bq)\partial_t(k,\bq))>\right)\\
&+&\partial_x\left(\frac{1}{k}<I(\d Q(k,\bq)\partial_x(k,\bq))>\right)\\
&-&\partial_t\left(<\KK(k,\bq)[B_T(k,\bq)\partial_t(k,\bq)]>\right)\\
&-&\partial_t\left(<\KK(k,\bq)[B_X(k,\bq)\partial_x(k,\bq)]>\right)\\
&-&\partial_x\left(c(k,\bq)<\KK(k,\bq)[B_T(k,\bq)\partial_t(k,\bq)]>\right)\\
&-&\partial_x\left(c(k,\bq)<\KK(k,\bq)[B_X(k,\bq)\partial_x(k,\bq)]>\right).
\end{array}
\right.
\end{equation}
In the next section, we will connect the issue of the well-posedness of system \eqref{whitham2_xt} with some second-order properties of spectral stabilty of periodic travelling-wave solutions to system \eqref{sv} under low-frequency perturbations. However we postpone nonlinear validation of this second-order modulation system to further work. Nevertheless in the last section we do explain what gain may be expected from its sudy.

\section{Spectral validation of modulation systems}

In this section, we carry out a spectral validation of systems \eqref{whitham_xt} and \eqref{whitham2_xt} by connecting their spectral properties to spectral properties of system \eqref{sv} in the low-frequency regime.

Let us fix $(\ks,\bqs)\in\R^\star\times\R$. Linearization will be performed for modulation systems  around $(\ks,\bqs)$, and correspondingly for the Saint-Venant equations around $(\Hs,\Qs)$. Moreover, as it is classical when studying the stability of a travelling wave, we will work in a co-moving frame, either $(x-\cs t,t)$ or $(\ks x+\oms t,t)$.

\subsection{The Whitham's system and Evans function}\label{subsec_Evans}

In this subsection, we show that the dispersion relation that determines the hyperbolicty of the Whitham's equations provides the principal part of the expansion of the Evans function, associated to the spectral stability of roll-waves. Namely we prove the following lemma.

\begin{lemma}\label{evansserre}
There is a $\Gamma\neq0$ such that 
$$
E(\lambda,e^{\nu})\ \stackrel{(0,0)}{=}\ \Gamma D(\lambda,\nu)\ +\ \mathcal{O}\left((|\lambda|+|\nu|)^3\right)
$$
where $E$ is the Evans function associated to the Saint-Venant equations (see \eqref{defEVsv} below) and $D$ is the dispersion of the Whitham's system (see \eqref{disp} below).
\end{lemma}

For that purpose, we follow the computation carried in \cite{N1} and consider a parametrization of viscous roll-waves by $(c,\bq)$ rather than $(k,\bq)$. The equivalence between the two parametrizations relies on assumption \eqref{cparam}. This assumption is fulfilled on the small-amplitude regime but degenerates in the small-viscosity limit. We will work right here in regimes where the assumption holds ; however note that some information may be obtained directly in $(k,\bq)$-variables as we will show with the Bloch analysis of the last subsection of this section.

In this subsection let us then denote $L(c,\bq)\in\R^\star$ a period and $H(\,\cdot\,;c,\bq)$ a $L(k,\bq)$-periodic function such that $H(L\,\cdot\,;c,\bq))$ is a $1$-periodic solution to \eqref{profil} with $k=[L(c,\bq)]^{-1}$. Accordingly $M(c,\bq)=\frac{1}{L(c,\bq)}\int_0^{L(c,\bq)} H(y;c,\bq)\d y$ and more generally 
\begin{equation}
<f>_L\ =\ \frac1L \int_0^L f\ .
\end{equation}

Since we are working with continuous roll-waves, we may choose a para\-metrization of $H$ such that $H(0;c,\bq)=M(c,\bq)$. Recall that this does not change the Whitham's system. This sligthly simplifies computations since then
$$
\d M\ =\ \left[-\frac{M}{L}+\frac{H(L)}{L}\right]\d L\,+\,<\d H>_L\ =\ <\d H>_L\ .
$$

Writting \eqref{whitham_xt} in frame $(x-\cs t,t)$ and linearizing yields
\begin{equation}\label{whit_ref}
\begin{array}{rclcl}
\displaystyle
\d\Ls\,\partial_t(c,\bq)&-&\Ls\partial_x c&=&0,\\
\displaystyle

<\d\Hs>_{\Ls}\,\partial_t(c,\bq)&+&\Ms\partial_x c-\partial_x \bq&=&0.
\end{array}
\end{equation}
or equivalently
\begin{equation}
\begin{array}{ccccl}
\displaystyle
\d\Ls\,\partial_t(c,\bq)&-&\Ls\partial_x c&=&0,\\
\displaystyle
\left[<\d\Hs>_{\Ls}+\frac{\Ms}{\Ls}\d\Ls\right]\,\partial_t(c,\bq)&-&\partial_x \bq&=&0.
\end{array}
\end{equation}
\noindent
Searching for solutions in the form $e^{\lambda t+\ks\nu x}(c_1,\bq_1)$, with $(c_1,\bq_1)$ constant and non-zero and $(\lambda,\nu)\in\C\times\txt{i}\R$, yields dispersion relation $0=D(\lambda,\nu)$ where
\begin{equation}\label{disp}
\displaystyle
\begin{array}{rcl}
D(\lambda,\nu)&:=&\lambda^2\left(\partial_c \Ls<\partial_{\bq}\Hs>_{\Ls}-\partial_{\bq} \Ls<\partial_c\Hs>_{\Ls}\right)\\
&&-\lambda\nu\left(\frac{\partial_c \Ls}{\Ls}+<\partial_{\bq} \Hs>_{\Ls}+\frac{\Ms\partial_{\bq} \Ls}{\Ls}\right)+\nu^2\ .
\end{array}
\end{equation}

Now we recall the construction of the Evans function given in \cite{N1}. Writting the Saint-Venant equations in the frame  $(x-\cs t,t)$ and linearizing gives a linear equation $\displaystyle \partial_t u-\As u=0$ for $u=(h,q)$, whose operator $\As$ has a spectrum on $L^\infty(\R)$ composed of $\lambda\in\mathbb{C}$ such that there exist functions $(h,q)$ and $\sigma\in\mathbb{S}^1$ such that
\begin{equation}\label{spec}\left\{
\begin{array}{rcl}
\displaystyle
\partial_x(q-\cs h)+\lambda h&=&0,\\
\displaystyle
\partial_x\left[\d G(\Hs,\Qs)(h,q)-\cs q\right]+\lambda q&=&\d S(\Hs,\Qs)(h,q)+\delta \partial_x^2q,
\end{array}\right.
\end{equation}
and
\begin{equation}
(q,h,h')(\Ls)\ =\ \sigma\ (q,h,h')(0)\ .
\end{equation}
\noindent
Setting $Y=(q,h,h')$, system \eqref{spec} can be written as a first order differential system with periodic coefficients:
\begin{equation}\label{edo1}
Y'=A(\lambda)Y.
\end{equation}
Let $\Psi(\,\cdot\,;\lambda)$ denote the fundamental resolvent matrix of \eqref{edo1}. Then the Evans function is
\begin{equation}\label{defEVsv}
E(\lambda,\sigma)=\det(\Psi(\Ls;\lambda)-\sigma \mathbb{Id}_{3}).
\end{equation}

To compute an expansion of $E$ in the neighbourhood of $(\lambda,\sigma)=(0, 1)$, we first choose a basis of solution to \eqref{edo1} for $\lambda=0$ and then continue this basis of solutions for small $\lambda$. We choose $Y_1(\,\cdot\,;0)=(\cs \Hs',\Hs',\Hs'')$ and $Y_3(\,\cdot\,;0)=(\cs \partial_{\bq}\Hs-1,\partial_{\bq}\Hs,\partial_{\bq}\Hs')$. Note that $Y_1$ and $Y_2$ are indeed independent. For $Y_2(\,\cdot\,;0)$, any choice completing the basis would do. But in order to compare with \cite{N1} let us choose  $Y_2(\,\cdot\,;0)=(\cs h_2,h_2,h'_2)$ where $h_2(\,\cdot\,;0)$ is a function such that $h_2(\Ls\,\cdot\,;0)$ belongs to the kernel of $\mathcal{L}_\star$ on $L^\infty(\R)$, provided by the Floquet analysis of profile equation \eqref{profil}, being associated to the Floquet multiplyier $\rho<1$ given by
$$
\displaystyle
\rho=\exp\left(\frac{1}{\delta \cs}\int_0^{\Ls}\d G(\Hs,\Qs)(1,\cs)-\cs^2\right).
$$ 
\noindent 
These eigenvectors can be continued analytically with respect to $\lambda$.

Let us denote, for any function f, $[f]_{\sigma}=f(\Ls)-\sigma f(0)$ and $[f]=[f]_1$, and perform the expansion
$$
\displaystyle
(q_j(\lambda),h_j(\lambda))\ =\ \sum_{l\in\N} \lambda^l (q_j^l,h_j^l)\ .
$$
Integrating \eqref{spec}$_{1}$ and performing the line substitution $L_1-\cs L_2\to L_1$ gives
{\setlength\arraycolsep{1pt}
\begin{eqnarray*}
\Delta^0(\lambda)E(\lambda,\sigma)&=&\left|\begin{array}{c}
                                \displaystyle                            
-\lambda\int_0^{\Ls}h_j(\lambda)+(\sigma-1)\bq_j(\lambda)\\
                                \displaystyle
[h_j(\lambda)]_\sigma\\
                                \displaystyle
[h_j(\lambda)']_\sigma 
                                \end{array}\right|_{1\leq j\leq3}\nonumber
\end{eqnarray*}}
with $\Delta^0(\lambda)=\det(Y_1(0;\lambda), Y_2(0;\lambda), Y_3(0;\lambda))$ and $\bq_j(\lambda)$ some constants such that 
$$
\cs h_j-q_j\ =\ \lambda I(h_j)+\bq_j\ .
$$
Note that $\bq_1(0)=\bq_2(0)=0$ and $\bq_3(0)=1$. Now, expanding $\bq_j(\lambda)=\sum_{l\in\N} \lambda^l \bq_j^l$ and looking at \eqref{spec}$_{2}$, one finds
\begin{equation*}
\LLs [h_1^1(\Ls\,\cdot\,;0)]\ =\ -\LLs [\partial_c\Hs(\Ls\,\cdot\,)]\ -\ (\bq_1^1-\Hs(0))\LLs[\partial_{\bq}\Hs(\Ls\,\cdot\,)] .
\end{equation*}
Therefore there are $\alpha$, $\beta$ such that
\begin{equation*}
h_1^1\ =\ -\partial_c\Hs\ -\ (\bq_1^1-\Hs(0))\partial_{\bq}\Hs\ +\ \alpha h_1^0\ +\ \beta h_2^0\ .
\end{equation*}

As a consequence, after column substitution $C_1-\lambda\beta C_2-\lambda(\bq_1^1-\Hs(0))C_3\to C_1$, comes
{\setlength\arraycolsep{2pt}
\begin{eqnarray*}
\Delta^0(0)E(\lambda,\sigma)&=&\left|\begin{array}{ccc}
                                \displaystyle                            
\lambda^2\int_0^{\Ls}\partial_c \Hs+\lambda(\sigma-1)\Hs(0) &-\lambda\int_0^{\Ls} h_2^0 &-\lambda\int_0^{\Ls}\partial_{\bq} \Hs+\sigma-1\\
                                \displaystyle
-\lambda\left[\partial_c \Hs\right]-(\sigma-1)h_1^0(0) &[h_2^0]& \left[\partial_{\bq} \Hs\right]\\
                                \displaystyle
-\lambda\left[\partial_c \Hs'\right]-(\sigma-1)h_1^0{}' (0) &[h_2^0{}']& \left[\partial_{\bq} \Hs'\right]\\
                                \end{array}\right|\nonumber\\
&+&\mathcal{O}\Big((|\lambda|+|\sigma-1|)^3\Big).
\end{eqnarray*}}
Differentiating $H(L(c,\bq),c,\bq)=H(0,c,\bq)$ and $H'(L(c,\bq),c,\bq)=H'(0,c,\bq)$ with respect to $(c,\bq)$ leads to
$$
\begin{array}{lll}
\displaystyle
[\partial_c \Hs]=-\partial_c \Ls\,\Hs'(0),\quad [\partial_c \Hs']=-\partial_c \Ls\,\Hs''(0),\\
\\
\displaystyle
[\partial_{\bq} \Hs]=-\partial_{\bq} \Ls\,\Hs'(0),\quad [\partial_{\bq} \Hs']=-\partial_{\bq} \Ls\,\Hs''(0),
\end{array}
$$
(whereas $[h_2^0]=(\rho-1)h_2(0;0)$ and  $[h_2^0{}']=(\rho-1)h_2{}'(0;0)$). This is now a straightforward computation to show that for some $\Gamma\neq0$
$$
\displaystyle
E(\lambda,e^{\nu})\ =\ \Gamma D(\lambda,\nu)\ +\ \mathcal{O}\left((|\lambda|+|\nu|)^3\right).
$$

\noindent
{\bf Remark:} the expansion of the Evans function found here is slightly different from the one derived in \cite{N1}: this latter one is not correct because of an uncorrect expansion of $[q_1(\lambda)-\cs h_1(\lambda)]+(1-\sigma)(q_1-\cs h_1)(0;\lambda)$ with respect to $\lambda$.

For the sake of completeness, in an appendix, we provide a more geometric description of this subsection that makes it more comparable with \cite{Serre}.

\subsection{Formal asymptotics in the Whitham's system}

To make some direct formal use of Lemma~\ref{evansserre}, in this subsection we study formally the hyperbolicity of \eqref{whitham_xt} in the regime of small viscosity, $\delta\to 0$, that is when viscous roll-waves are close to Dressler roll-waves. For that purpose, we use $(h_+,\bq)$-para\-metrization (as in Proposition~\ref{Dressler}) and will suppose that all concerned quantities are regular with respect to $(\delta,h_+,\bq)$. Moreover, we write equations in the reference co-moving frame $(x-\cs t,t)$. One obtains asymptotically
$$
\begin{array}{rcccccl}
\displaystyle
\frac{H_c(\bqs)}{F^2}\mathcal{P}'_\star\ \partial_t h_+ 
&+&\displaystyle \frac{H_c'(\bqs)}{F^2}\mathcal{P}_\star\ \partial_t \bq
&-&\Ls c^*{}'(\bqs)\ \partial_x \bq&=&0,\\
\displaystyle
H_c(\bqs)\left(\frac{\mathcal{Q}}{\mathcal{F}}\right)'_\star\ \partial_t h_+ 
&+& \displaystyle H_c'(\bqs)\left(\frac{\mathcal{Q}}{\mathcal{F}}\right)_\star\ \partial_t \bq
&+&\left(\Ms c^*{}'(\bqs)-1\right)\ \partial_x \bq&=&0.
\end{array}
$$
so that $(h_+,\bq)(x)=e^{\nu x}(1,0)$ defines a solution for any $\nu$. Therefore $0$ is always a trivial eigenvalue and the $2\times2$-system is always hyperbolic. The Whitham's system being {\it always} hyperbolic in the vanishing viscosity limit $\delta\to0$, one can expect, using a pertubation argument, that the Whitham's equations are hyperbolic for $\delta>0$ sufficiently small. This is in contrast with the results of Boudlal and Liapidevskii \cite{Boudlal} that suggested that inviscid roll-waves are stable under long wavelength perturbation (in the sense that the Whitham's system is hyperbolic) only for roll-waves of limited periods. This discrepancy may come from the fact that the modulation procedure was not carried out properly in \cite{Boudlal}. Besides, the formal asymptotics carried out here confirms \emph{formally} the spectral analysis of inviscid roll-waves performed rigorously in \cite{N2}. 

Let us emphasize that trivial hyperbolicity in the formal asymptotic system comes from the fact that wave speed $c$ does not depend on $k$ (nor on $h_+$) but only on $\bq$.

In contrast, a discussion in the neighbourhood of Hopf's bifurcation points using a perturbation argument from the constant case would show instability (see \cite{Ba_Jo_Ro_Zu}). Yet, some numerical evidence shows that somewhere between bifurcation points and limiting homiclinics there are some stable roll-waves and that, in the weak stable sense of hyperbolicity of the Whitham's system, roll-waves are even stable up to the homoclinic travelling-waves. A detailed discussion of these former points (instablity close to bifurcation, numerical check of stability for roll-waves) may be found in \cite{Ba_Jo_No_Ro_Zu_2, Ba_Jo_No_Ro_Zu_1,Ba_Jo_Ro_Zu}.

\subsection{Second-order modulation and Bloch-wave analysis}

In this subsection, we validate spectral properties of \eqref{whitham2_xt} by a spectral Bloch-wave analysis of \eqref{sv}. Also we obtain relations between eigenvectors of the Saint-Venant system and those of the modulation systems (both \eqref{whitham_xt} and \eqref{whitham2_xt}).

We are thus lead to perform a Bloch study of the operator $\As$, defined as
\begin{equation}\label{As}
\displaystyle
\As (h,q)=\left(\begin{array}{c}\displaystyle\ks\partial_y(\cs h-q)\\
\displaystyle\partial_h\GGs[h]+\partial_q\GGs[q]-\oms\partial_y q\end{array}\right)\ .
\end{equation}
Yet, the reader is referred to next section for a definition of the Bloch transform. Let us only introduce, for $l\in[-\pi, \pi]$, the operator $\Asc(l)$ acting on $L^2_{\textrm{per}}(\R)$ and defined by
$$
\displaystyle
[\Asc(l)(f)](y)=e^{-il y}[\As (e^{il \cdot}f(\cdot))](y)\ .
$$

The following lemma is a version of Lemma~2.1 in \cite{Jo_Zu_No}, stated for the Saint-Venant equations in Lagrangian coordinates and used to prove that, for viscous roll-waves, linear stability implies nonlinear stability.

\begin{lemma}\label{eigenvect0}
The critical eigenvalues $\lambda_j, j=1,2$ of $\Asc$ are analytic functions of the Floquet number $l$. 

The Jordan structure of the zero eigenspace of $\Asc(0)=\As$ consists of a $1$-dimensional kernel and a single Jordan chain of height $2$. The left kernel of $\As$ is spanned by the constant function $(1, 0)$ and $(\Hs',\Qs')$ spans the right eigendirection lying at the base of the Jordan chain.

Moreover, for $|l|$ sufficiently small, there exist dual right and left eigenfunctions $w_j(\cdot,l)$ and $\tilde{w}_j(\cdot,l)$ of $\Asc(l)$ associated with $\lambda_j(l)$, for $j=1,2$, of form 
$$
w_j=\sum_{k=1}^2 \beta_{j,k}v_k\ ,\quad\tilde{w}_j=\sum_{k=1}^2\tilde{\beta}_{j,k}\tilde{v}_k
$$ 
where\begin{itemize} 
\item $(v_k)_{k=1,2}$ and $(\tilde{v}_k)_{k=1,2}$ are dual bases of the total eigenspace of $\Asc(l)$ associated with small eigenvalues, analytic in $l$, and such that
$$\displaystyle
\tilde{v}_2(\cdot,0)=(1,0)\ ,\quad v_1(\cdot,0)=(\Hs',\Qs')\ ;
$$
\item $(l^{-1}\tilde{\beta}_{j,1}(l),\tilde{\beta}_{j,2}(l))$, $j=1,2$, and $(l\beta_{j,1}(l),\beta_{j,2}(l))$, $j=1,2$ are analytic in $l$.
\end{itemize}
\end{lemma}

The role of $(1,0)$ is a direct consequence of the fact that the first equation of the Saint-Venant system is a conservation law. As already pointed out, the role of $(\Hs',\Qs')$ stems form translational invariance of system \eqref{sv}. The rest of the lemma, the analyticity issue, may be obtained in a standard way (see \cite{Jo_Zu_No}) and is directly related to the existence of an averaged modulation system.

We come to the main part of this subsection.

\begin{lemma}\label{eigenvect1}
For $j=1,2$,
$$
\displaystyle
w_j(\cdot,l)=
\frac{1}{il}\left(k^0_j+ilk^1_j\right)
\left(\begin{array}{c}
\displaystyle
\frac{1}{\ks}\Hs'+il\partial_k\Hs\\
\displaystyle
\frac{1}{\ks}\Qs'+il\partial_k\Qs
\end{array}\right)
+\bq^0_j
\left(\begin{array}{c}
\displaystyle
\partial_{\bq}\Hs\\
\displaystyle
\partial_{\bq}\Qs
\end{array}\right)
+\mathcal{O}(l)
$$
with $(k^0_j,\bq^0_j)$ the first term in a low-frequency expansion of a corresponding eigenvector of the linearized Whitham's system. Thereby $v_1$ and $v_2$ may be choosed so that
$$
\displaystyle
v_1(\cdot,l)=
\left(\begin{array}{c}
\displaystyle
\frac{1}{\ks}\Hs'+il\partial_k\Hs\\  
\displaystyle
\frac{1}{\ks}\Qs'+il\partial_k\Qs
\end{array}\right)
+\mathcal{O}(l^2),
\quad 
v_2(\cdot,l)=
\frac{1}{\partial_{\bq}\Ms}
\left(\begin{array}{c}
\displaystyle
\partial_{\bq}\Hs\\
\displaystyle
\partial_{\bq}\Qs
\end{array}\right)
+\mathcal{O}(l).
$$
\end{lemma}

\begin{lemma}\label{eigenvect2}
For $j=1,2$,
$$
\lambda_j (l)\ =\ \mu_j(l)\ +\ \mathcal{O}(l^3)
$$
where $\mu_j(l)$ is a Fourier eigenvalue of system~\eqref{whitham2_xt} corresponding to frequency~$l$.
\end{lemma}

Lemma~\ref{eigenvect1} is a consequence of the proof of Lemma~\ref{eigenvect2}, so we focus on the proof of the latter. Moreover, we choose to write this proof in a wave analysis spirit. Lemma~\ref{eigenvect2} also comes with a better description of critical eigenvectors but we do not write it here.

Writting \eqref{whitham2_xt} in frame $(\ks x+\oms t,t)$ and linearizing yields, after some simplification in the second equation with help of the first one,
\begin{equation}\label{whitham2_linear}
\left\{
\begin{array}{rcl}
\partial_t k+\ks^2\d\cs\partial_x (k,\bq)
&=&-\ \ks\partial_x\left(\ps[\BTs\partial_t(k,\bq)]\right)\\
&-&\ks\partial_x\left(\ps[\ks\BBXs\partial_x(k,\bq)]\right)\\
\d\Ms\partial_t (k,\bq)
-\frac{\Ms}{\ks}\partial_t k \\
=\ \ks\partial_x\left(\frac{1}{\ks}<I(\d \Hs\partial_t(k,\bq))>\right)
&+&\ks\partial_x\left(\frac{1}{\ks}<I(\ks\Hs)>\d \cs\partial_x(k,\bq)\right)\\
&-&\ks\partial_x\left(\frac{1}{\ks}<\ks I(\un)>\partial_x\bq\right)\\
-\ \partial_t\left(<\KKs[\BTs\partial_t(k,\bq)]>\right)
&-&\partial_t\left(<\KKs[\ks\BBXs\partial_x(k,\bq)]>\right)
\end{array}
\right.
\end{equation}
with
\begin{equation}
\BBXs\ =\ \BXs\ -\ \cs\BTs\ .
\end{equation}

Looking for solutions of type $e^{\lambda t+\nu x}(k(\nu),\bq(\nu))$ with constant $(k(\nu),\bq(\nu))$ expanded into $\sum_{l\in\N}\nu^l (k^l,\bq^l)$, one finds at order $\mathcal{O}((|\lambda|+|\nu|)^2)$
\begin{equation*}
\lambda \left[\begin{array}{cc}1&0\\ \partial_k\Ms-\frac{\Ms}{\ks}&\partial_{\bq}\Ms\end{array}\right]\left(\begin{array}{c}k^0\\\bq^0\end{array}\right)
\ +\ \ks\nu \left[\begin{array}{cc}\ks\partial_k\cs&\ks\partial_{\bq}\cs\\ 0&-1\end{array}\right]
\left(\begin{array}{c}k^0\\\bq^0\end{array}\right)\ =\ \left(\begin{array}{c}0\\0\end{array}\right).
\end{equation*}
From \eqref{evolution} stems that $\lambda=0$ is the unique eigenvalue corresponding to $\nu=0$. We now write
$$
\displaystyle
\lambda\ =\ \ks\nu\lambda(\nu),\quad \lambda(\nu)\ =\ \sum_{l\in\N}\nu^l\lambda^l
$$
to get further information.

With this notations, we obtain
\begin{equation*}
\lambda^0 \left[\begin{array}{cc}1&0\\ \partial_k\Ms-\frac{\Ms}{\ks}&\partial_{\bq}\Ms\end{array}\right]\left(\begin{array}{c}k^0\\\bq^0\end{array}\right)
\ +\ \left[\begin{array}{cc}\ks\partial_k\cs&\ks\partial_{\bq}\cs\\ 0&-1\end{array}\right]
\left(\begin{array}{c}k^0\\\bq^0\end{array}\right)\ =\ \left(\begin{array}{c}0\\0\end{array}\right).
\end{equation*}
Naturally we recover the same dipersion relation as for \eqref{whitham_xt}, thus spectral Bloch-wave validation of \eqref{whitham2_xt} contains some spectral validation of \eqref{whitham_xt}.

At next order, holds
\begin{equation*}\label{whitham2_spec}
\left\{
\begin{array}{l}
\lambda^1 k^0+\lambda^0 k^1+\ks\d\cs(k^1,\bq^1)\\
=-\ks\ps[\left(\lambda^0\BTs+\BBXs\right)(k^0,\bq^0)]\\
\lambda^1\d\Ms(k^0,\bq^0)+\lambda^0\d\Ms(k^1,\bq^1)
-\lambda^1\frac{\Ms}{\ks}k^0-\lambda^0\frac{\Ms}{\ks}k^1\\
=\ \ks\left(\frac{\lambda^0}{\ks}<I(\d \Hs(k^0,\bq^0))>\right)
-\ks\left(\frac{\lambda^0}{\ks}<I(\Hs)>\frac{k^0}{\ks}\right)\\
-\ks\left(\frac{1}{\ks}<I(\un)>\bq^0\right)
-\ \ks\lambda^0\left(<\KKs[\left(\lambda^0\BTs+\BBXs\right)(k^0,\bq^0)]>\right)
\end{array}
\right.
\end{equation*}
since
\begin{equation}\label{simply}
\lambda^0 k^0\ +\ \ks\d\cs (k^0,\bq^0)\ =\ 0\ .
\end{equation}
Moreover, using again \eqref{simply}, one may obtain
\begin{equation*}\label{Blambda}
\begin{array}{l}
\left(\lambda^0\BTs+\BBXs\right)(k^0,\bq^0)\\
=\ \left[-\cs^2\d\Hs-\d G(\d\Hs,\cs\d\Hs)-2\delta\ks\d\Hs'\right](k^0,\bq^0)\\
-\ \delta\cs\Hs' k^0\ -\ \left[\frac{1}{\ks}\partial_q\GG(\Hs,\Qs;\ks)[I(\un)]+\partial_q G(\Hs,\Qs)\right]\bq_0\\
-\ \lambda^0\frac{k^0}{\ks}\left[\frac{1}{\ks}\partial_q\GG(\Hs,\Qs;\ks)[I(\Hs)]+\partial_q G(\Hs,\Qs)\Hs-2\delta\ks\Hs'\right]\\
+\ \lambda^0\left[2\cs\d\Hs+\frac{1}{\ks}\partial_q\GG(\Hs,\Qs;\ks)[I(\d\Hs)]\right](k^0,\bq^0)\\
-\ \lambda^0 \bq^0\ -\ (\lambda^0)^2\frac{k^0}{\ks}\Hs\\
=\ \Bs^0(k^0,\bq^0)\ +\ \lambda^0\Bs^1(k^0,\bq^0)\ +\ (\lambda^0)^2\Bs^2(k^0,\bq^0)\ 
=\ \Bs[\lambda^0](k^0,\bq^0)\ .
\end{array}
\end{equation*}

We now carry out a spectral Bloch-wave analysis of system \eqref{sv}. Writting the system in frame $(\ks x+\oms t,t)$ and linearizing gives
\begin{equation}\label{Bloch}
\left\{\begin{array}{rcl}
\partial_t h&=&\ks\partial_x(\cs h-q),\\
\partial_t q&=&-\oms \partial_x q\ +\ \partial_h\GG(\Hs,\Qs;\ks)[h]\ +\ \partial_q\GG(\Hs,\Qs;\ks)[q]
\end{array}\right.
\end{equation}
that is
$$
\partial_t (h,q)\ =\ \As (h,q)\ .
$$
We look for solutions to \eqref{Bloch} in the form $e^{\lambda t+\nu x}(h_\nu(x),q_\nu (x))$ with $(h_\nu,q_\nu)$ $1$-periodic functions. We are only interested in spectrum near $(\lambda,\nu)=(0,0)$. Therefore we set 
$$
\displaystyle
\lambda\ =\ \ks\nu\lambda(\nu),\quad \lambda(\nu)\ =\ \sum_{l\in\N}\nu^l\lambda^l
$$
and 
$$
\displaystyle
(h_\nu,q_\nu)\ =\ \sum_{l\in\N}\nu^l(h_l,q_l)\ .
$$

First, there must be $\bq\in\R$ such that $\cs h_0-q_0=\bq_0$ and then 
$$
\LLs h_0\ =\ \bq\Aqs\ =\ \LLs\partial_{\bq}\Hs\ - \bq \partial_{\bq}\oms \Aoms
$$
which in turn imposes $\bq=0$ (since $\partial_{\bq}\oms\neq0$). Hence there is $k^0\in\R^\star$ such that
$$
h^0\ =\ \frac{k^0}{\ks} \Hs',\quad q^0\ =\ \cs h^0\ =\ \frac{k^0}{\ks} \Hs'\ .
$$

Now there is $\bq^0\in\R$ such that
$$
\cs h^1-q^1\ =\ \bq^0\ +\ \lambda^0\,\frac{k^0}{\ks}\Hs
$$
and then
$$
\LLs h^1\ =\ k^0\AAks\ +\ \bq^0\Aqs\ +\ \lambda^0\,k^0\Aoms\ .
$$
This forces
\begin{equation}\label{Bloch_k0}
\lambda^0\,k^0\ +\ \ks\d\cs(k^0,\bq^0)\ =\ 0
\end{equation}
and the existence of $k^1\in\R$ such that
$$
h^1\ =\ \d\Hs(k^0,\bq^0)\ +\ \frac{k^1}{\ks}\Hs'\ .
$$

Then equation
$$
\partial_x(\cs h^2-q^2)\ =\ -(\cs h^1-q^1)\ +\ \lambda^1 h^0\ +\ \lambda^0 h^1
$$
implies
$$
\cs <h^1>-<q^1>\ =\ \lambda^0 <h^1>
$$
thus
\begin{equation}\label{Bloch_q0}
\lambda^0\left(\d\Ms(k^0,\bq^0)-\frac{\Ms}{\ks}k^0\right)\ -\ \bq^0\ =\ 0\ .
\end{equation}
Note that equations (\ref{Bloch_k0},\ref{Bloch_q0}) already provides a first-order spectral justification. Moreover there exists $\bq^1\in\R$ such that
$$
\begin{array}{rcl}
\cs h^2-q^2&=&\bq^1\ -\ \bq^0 I(\un)\ -\ \lambda^0\frac{k^0}{\ks}I(\Hs)\\
&+&\lambda^1\,\frac{k^0}{\ks}\Hs\ +\ \lambda^0\,\frac{k^1}{\ks}\Hs\ +\ \lambda^0\,I(\d\Hs(k^0,\bq^0))\ .
\end{array}
$$
Now
$$
\LLs h^2\ =\ k^1\AAks\ +\ \bq^1\Aqs\ +\ (\lambda^1\,k^0+\lambda^0\,k^1)\Aoms\ +\ R^1
$$
with
$$
R^1\ =\ \ks \Bs[\lambda^0](k^0,\bq^0)\ .
$$
This sets 
\begin{equation}\label{Bloch_k1}
\lambda^1\,k^0\ +\ \lambda^0\,k^1\ +\ \ks\d\cs(k^0,\bq^0)\ =\ -\ks\ps[\Bs[\lambda^0](k^0,\bq^0)]\ .
\end{equation}
Moreover there exists $k^2$ such that
$$
h^2\ =\ \d\Hs(k^1,\bq^1)\ +\ \KKs R^1\ +\ \frac{k^2}{\ks}\Hs'\ .
$$
Now equation
$$
\cs <h^2>-<q^2>\ =\ \lambda^0 <h^2>\ +\ \lambda^1 <h^1>
$$
leads to
\begin{equation}\label{Bloch_q1}
\begin{array}{l}
\lambda^1\d\Ms(k^0,\bq^0)\ +\ \lambda^0\d\Ms(k^1,\bq^1)\ +\ \lambda^0<\KKs \Bs[\lambda^0](k^0,\bq^0)>\\
=\ \lambda^1\frac{\Ms}{\ks}k^0\ +\ \lambda^0\frac{\Ms}{\ks}k^1\ +\ \bq^1\ +\ \lambda^0<I(\d\Hs(k^0,\bq^0))>\\
-\ \lambda^0<I(\Hs)>\frac{k^0}{\ks}\ -\ \bq^0<I(\un)>
\end{array}
\end{equation}
which completes our spectral justification.

Indeed writting a linear system for $(k^0,\bq^0)+\nu(k^1,\bq^1)$, which is non-zero, leads to the same dispersion relation for $\lambda^0+\nu\lambda^1$ in both cases.

Let us comment somewhat on the spectrum we just described. For $\nu=0$ in the Bloch-wave analysis, $\lambda=0$ is an eigenvalue corresponding to a $2\times2$-Jordan block, with $(\Hs',\Qs')$ as an eigenvector (and $\frac{1}{\ks}(\partial_c\Hs,\partial_c\Qs)$ as its antecedent in the $(c,\bq)$-parametrization). For $\nu\in \txt{i}\R^\star$ small, two eigenvalues emerge from $0$, tangent to the imaginary axis when hyperbolicity of \eqref{whitham_xt} is met, with first order $\ks\nu\lambda^0$, $\lambda^0$ being well-described by both modulation systems. Curvatures of the eigenvalue curves may then be extracted from system \eqref{whitham2_xt}.

\section{Nonlinear validation\\ of the Whitham's system}

In this section, we prove the existence of a family of solutions to the shallow water equations, close to a given roll-wave, converging towards a modulated roll-wave profile described at first order by a solution to the {\it inviscid} Whitham's system.

Firstly we enumerate here all our assumptions. Not all of them will be recalled in Theorem~\ref{main}. Beyond the assumptions of smooth $(k,\bq)$-parametrization and that the Whitham's system is of evolution-type (see \eqref{evolution}) and hyperbolic, we will also ask for :
\begin{itemize}
\item the spectrum of $\As$ is of upper bounded real part ;
\item for any $l\in [-\pi,\pi]$ and any cut-off parameter $l_1$, the non critical part $\Asc(l)\PsFS$ of $\Asc(l)$ is invertible (see \eqref{Qcc}, \eqref{cut-off} and \eqref{FSMF} for definitions).
\end{itemize}

We write the latter points as assumptions to emphasize what we really use in our proof. But these former assumptions may be removed. For instance, the upper-boundedness is a consequence of high-frequency estimates in \cite{N1}. Alternatively, invertibility assumption, as hyperbolicity, can be rigourously reduced to numerical investigation in an explicit finite box of eigenvalue phase space and then numerically checked with techniques in \cite{Ba_Jo_No_Ro_Zu_2, Ba_Jo_No_Ro_Zu_1,Ba_Jo_Ro_Zu}.

\subsection{Spaces}

Following the strategy introduced in \cite{Me_Sch} for the Ginzburg-Landau equations and developped in \cite{D3S} for reaction-diffusion systems, we will prove the convergence to a roll-wave profile in a set of \emph{analytic} functions. Indeed, in the hyperbolic scaling considered here, there is no smoothing effect arising from equations whereas, in the modulation context, some terms are neglected precisely because they contain more derivatives. We will correspondingly restrict the class of admissible solutions to the Whitham's system and to the Saint-Venant equations.

Let $a>0$ and $m\in\R^+$. We first introduce a space for solutions to the Saint-Venant system. Let us define
$$
\displaystyle
L_{\mathcal{J}}(a,m)=\left\{v\in L^1([-\pi,\pi], H^m_{\txt{per}})\middle|\:\int_{-\pi}^{\pi}\|v(\cdot,l)\|_{H^m_{\txt{per}}}e^{a|l|}dl<\infty\right\}
$$
where $H^m_{\txt{per}}$ denotes the classical Sobolev space of $1$-periodic functions. Also, for any Schwartz-class function $u$, $\check{u}$ is the {\it Bloch transform} of $u$ defined by, for any $(y,l)\in\R\times[-\pi, \pi]$
$$
\displaystyle
\check{u}(y,l)=\mathcal{J}u(y,l)=\sum_{j\in\Z}e^{i2\pi jy}\mathcal{F}u(l+2\pi j),
$$
$\hat{u}=\mathcal{F}u$ being the Fourier transform of $u$, explicitely for $l\in\R$
$$
\hat{u}(l)\ =\ \mathcal{F}u(l)\ =\ \frac{1}{\sqrt{2\pi}}\int_\R e^{-ilx}u(x)\d x\ .
$$
Note that a justification for restricting attention to $\R\times[-\pi,\pi]$ is that extending definition to $\R\times\R$ would lead to : for any $(y,l)\in\R\times\R$,
$$
\check{u}(y,l+2\pi)\ =\ e^{-i2\pi y}\check{u}(y,l)\ .
$$
Moreover, the Bloch transform comes with an inverse formula
$$
u(x)\ =\ \frac{1}{\sqrt{2\pi}}\int_{-\pi}^\pi e^{ilx}\check{u}(x,l)\d l
$$
and a Plancherel formula
$$
\|u\|_{L^2(\R)}=\|\check{u}\|_{L^2([-\pi,\pi],L^2_{\txt{per}})}.
$$
Admissible solutions will be considered in the Banach space
$$
\displaystyle
\mathcal{X}^{a}_m=\left\{u\,:\,\R\to \C^2\middle|\:\check{u}\in L_{\mathcal{J}}(a,m)\right\}
$$
endowed with norm $\xnorm{\,\cdot\,}{a}{m}$ defined as, for any Schwartz $u$,
$$
\xnorm{u}{a}{m}=\|\check{u}\|_{L_{\mathcal{J}}(a,m)}=\int_{-\pi}^{\pi}\|\check{u}(\cdot,l)\|_{H^m_{\txt{per}}}e^{a|l|}dl\ .
$$
Due to Sobolev embedding theorems, $\mathcal{X}^{a}_m$ is an algebra when equipped with usual multiplication provided $m\geq 1$. Namely, if $m\geq1$, there is a $C(m)>0$ (independent of $a$) such that for any $u,v\in\mathcal{X}^{a}_m$ 
$$
\displaystyle \xnorm{uv}{a}{m}\leq C(m)\xnorm{u}{a}{m}\xnorm{v}{a}{m}\ .
$$ 

Correspondingly we introduce a space for solutions to the Whitham's system. First
$$
\displaystyle
L_{\mathcal{F}}(a,m)=\left\{v\in L^1(\R,\C)\,\middle|\:\int_{\R}|v(l)|(1+|l|)^m e^{a|l|}dl<\infty\right\}
$$
and 
$$
\displaystyle
\mathcal{Y}^{a}_m=\left\{u\,:\,\R\to \C\,\middle|\:\hat{u}\in L_{\mathcal{F}}(a,m)\right\}
$$
endowed with $\ynorm{\cdot}{a}{m}$ :
$$
\ynorm{u}{a}{m}=\|\hat{u}\|_{L_{\mathcal{F}}(a,m)}=\int_{\R}|\hat{u}(l)|(1+|l|)^m e^{a|l|}dl\ .
$$
Note that such $u$ are analytic on strip $\{z\in\C||\Im z|<a\}$. Moreover, when $m\geq1$, $\mathcal{Y}^{a}_m$ is also an algebra.

An important link between two kinds of settings is provided by the fact that if $u$ is $1$-periodic and the Fourier transform of $v$ is supported in $[-\pi,\pi]$ then
$$
\mathcal{J}(uv)(y,l)\ =\ u(y)\ \mathcal{F}(v)(l)\ 
$$
and in particular, for such a low-frequency $v$, $\mathcal{J}(v)(y,l)=\mathcal{F}(v)(l)$. Conversely, averaging in $y$ yields, for a general $u$,
$$
\left<\mathcal{J}(u)(\cdot,l)\right>\ =\ \mathcal{F}(u)(l)\ 
$$
for any $l\in\R$ (with an extended Bloch transform).

To make our theorem more readible we also introduce more common spaces : uniformly local Sobolev spaces. Let us first introduce an intermediate space
$$
\tilde{H}^m\ =\ \left\{u\in L^2_{\txt{loc}}\ \middle|\ \sup_{x\in\R}\|u_{|[x,x+1]}\|_{H^m([x,x+1])}<\infty\right\}
$$
endowed with norm
$$
\hulnorm{u}{m}\ =\ \sup_{x\in\R}\|u_{|[x,x+1]}\|_{H^m([x,x+1])}\ .
$$
Then we define the subspace
$$
H^m_{\txt{ul}}\ =\ \left\{u\in\tilde{H}^m\ \middle|\ \begin{array}{l}\R\to\tilde{H}^m,\\\tau\mapsto u(\cdot-\tau)\end{array}\ \txt{is continuous}\right\}\ .
$$

\subsection{Main statement}

Before stating the main result of this section, which we will prove in the following, we still need to introduce some change of variables. Let $(\ks,\bqs)\in\R^\star_+\times\R$ and then write equations in usual frame $(\ks x+\oms t,t)$. Then, for any phase $\varphi$, we introduce a time-dependent change of variable $\Xphi$ defined by
$$
\Xphi(y,t)\ =\ y\ -\ \varphi(y,t)\ .
$$
If $\partial_y\varphi$ is small enough, it can be inverted into $\Yphi$ satisfying
$$
\Xphi(\Yphi(x,t),t)\ =\ x\ .
$$
Note that
\begin{eqnarray}\label{Yphix}
\partial_x \Yphi(x,t)&=&\frac{1}{1-\partial_y\varphi(\Yphi(x,t),t)}\\ \label{Yphit}
\partial_t \Yphi(x,t)&=&\frac{\partial_t \varphi(\Yphi(x,t),t)}{1-\partial_y\varphi(\Yphi(x,t),t)}\ .
\end{eqnarray}
In particular, whenever $\partial_y\varphi$ is small, $\partial_x\Yphi(x,t)$ is close to $1+\partial_y\varphi(\Yphi(x,t),t)$, corresponding local wavenumber is close to $\ks(1+\partial_y\varphi(\Yphi(x,t),t))$, and corresponding local frequency close to 
$$
\partial_t\varphi(\Yphi(x,t),t)+\oms(1+\partial_y\varphi(\Yphi(x,t),t))\ .
$$
Moreover all derivatives are written as functions of $\Yphi(x,t)$ and $t$.

Now the following theorem provides us with a nonlinear justification of the Whitham's equations
\begin{equation}\label{just_whit}\left\{
\begin{array}{ll}
\displaystyle
\partial_T k+\ks\partial_X(k c(k,\bq)-\cs k)=0\\
\displaystyle
\partial_T M(k,\bq)+\ks\partial_X\big(N(k,\bq)-\cs M(k,\bq)\big)=0
\end{array}\right.
\end{equation}
here written in frame $(\ks X+\oms T,T)$.

\begin{theorem}\label{main}
Let $(\ks,\bqs)\in\R^\star_+\times\R$. Assume that in a neighbourhood of $(\ks,\bqs)$ system \eqref{just_whit} is of evolution type and strictly hyperbolic.\\
For any $a>0$, $m\geq 3$ and $M\geq 1$, there exist positive $\varepsilon_1,\eta_1, C_1$ and $T_1$ such that, for any $T_0\in]0,T_1]$, for any solution $(k,\bq)$ to \eqref{just_whit} on $[0, T_0]$ satisfying 
$$
\displaystyle
\sup_{T\in[0, T_0]}\ynorm{(k,\bq)(\cdot,T)-(\ks,\bqs)}{a}{0}\leq \eta_1
$$
and for all $\varepsilon\in]0,\varepsilon_1[$, there exist $\big((k_{\varepsilon},\bq_{\varepsilon}),r^s_{\varepsilon}\big)$ and $\varphi_{0,\varepsilon}$ such that
$$
\begin{array}{ll}
\displaystyle
\sup_{t\in [0, T_0/\varepsilon]}\hulnorm{(k_{\varepsilon},\bq_{\varepsilon})(\cdot,t)-(k,\bq)(\varepsilon \cdot,\varepsilon t)}{m}
\leq C_1\left[\varepsilon+\eta^2\right],\\
\displaystyle
\sup_{t\in [0, T_0/\varepsilon]}\hulnorm{r^s_{\varepsilon}(\cdot,t)}{m}\leq C_1\eta^2,\\
\displaystyle
\sup_{t\in[0,T_0/\varepsilon]}|\varphi_{0,\varepsilon}(t)|\leq C_1 \frac{\eta^2}{\varepsilon}
\end{array}
$$
where
$$
\eta\ =\ \sup_{T\in[0,T_0]}\hulnorm{(k,\bq)(\cdot,T)-(\ks,\bqs)}{m}
$$
and a solution $(h,q)$ to the Saint-Venant system \eqref{sv} such that
$$
\displaystyle
\sup_{t\in[0, T_0/\varepsilon]}\sup_{x\in\R}|(h,q)(x,t)-(H_{app},Q_{app})(Y^{\varphi_\varepsilon}(x,t),t)|\leq C_1\varepsilon^M,
$$
where 
$$
\displaystyle
(H_{app},Q_{app})(y,t)=(H,Q)(y,k_{\varepsilon}(y,t),\bq_{\varepsilon}(y,t))+r^s_{\varepsilon}(y,t),
$$
and $\varphi_\varepsilon(y,t)=\varphi_{0,\varepsilon}(t)+\int_0^y (\frac{k_{\varepsilon}(z,t)}{\ks}-1)\d z$.
\end{theorem}

Our result is partially stated in $\hulnorm{\cdot}{m}$ norm, whereas we will work with Fourier, Bloch, or mixed types, multipliers. To fill this gap we need a multiplier theorem and thus refer to Lemma~$5$ in \cite{Schneider_Error} (also stated in \cite{D3S} as Lemma~$3.6$).

\subsection{Separation of critical modes}

Our proof of the above theorem starts simultaneously rewriting the shallow water equations in an appropriate form so as to separate critical modes from others and plugging into the equations a \emph{roll-wave ansatz} that we will modulate afterwards when looking for long wavelength approximate solutions.

Recall that, in frame $(\ks x+\oms t,t)$, the Saint Venant system is written as
\begin{equation}\label{stv2w}
\left\{
\begin{array}{rcl}
\displaystyle
\partial_t h+\ks\partial_x (q-\cs h)&=&0\\
\displaystyle
\partial_t q+\ks\partial_x\big(G(h,q)-\cs q\big)&=&S(h,q)+\delta\ks^2\partial^2_x q
\end{array}\right.
\end{equation}
or equivalently
\begin{equation*}
\left\{
\begin{array}{rcl}
\displaystyle
\partial_t h+\ks\partial_x (q-\cs h)&=&0\\
\displaystyle
\partial_t q+\oms\partial_x q-\GG(h,q;\ks)&=&0
\end{array}\right.\ .
\end{equation*}

Following \cite{D3S}, we introduce a roll-wave ansatz for $(h,q)(x,t)$,
$$
\displaystyle 
(H,Q)\left(\Yphi(x,t);\frac{\ks}{1-\partial_y\varphi(\Yphi(x,t),t)},\bqs+\bq(\Yphi(x,t),t)\right)+(\hh,\qq)(\Yphi(x,t),t).
$$ 
Since $\partial_y\varphi$ should remain small, $(h,q)(\Xphi(y,t),t)$ is then well-approximated by
$$
(\Hs,\Qs)(y)+(\d\Hs,\d\Qs)(\ks\partial_y\varphi(y,t),\bq(y,t))(y)+(\hh,\qq)(y,t).
$$

Using (\ref{Yphix},\ref{Yphit}) in a chain rule differentiation turns \eqref{stv2w} into a set of equations for $(\hh,\qq)$ and $(\varphi,\bq)$. We added two more unknowns and therefore should later add two more constraints. These constraints will perform a separation of low-Floquet critical modes.

As an example, note that
$$
\GG\left((H,Q)\left(\Yphi(\cdot,t);\frac{\ks}{1-\partial_y\varphi(\Yphi(\cdot,t),t)},\bqs+\bq(\Yphi(\cdot,t),t)\right);\ks\right)(\Xphi(y,t),t)
$$
would turn into
$$
\displaystyle
\GGs+\partial_h\GGs[\d\Hs(\ks\partial_y\varphi,\bq)]+\partial_q\GGs[\d\Qs(\ks\partial_y\varphi,\bq)]+\partial_k\GGs\ks\partial_y\varphi+\delta\ks^2\Qs'\partial_y^2\varphi+\GG_{\mathcal{R}}
$$
(taken in (y,t)), with $\GG_{\mathcal{R}}$ at least quadratic in $(\ks\partial_y \varphi,\bq)$.

Setting $u=(\hh,\qq)$, the shallow water system \eqref{stv2w} leads to
\begin{equation}\label{stv2w1}
\begin{array}{rclcl}\displaystyle
[\BWT_0+\BWT_1(u,\ks\partial_y\phi,\bq)](\ks\partial_t \varphi,\partial_t\bq)&+&\partial_t u&&\\
\displaystyle
-\BWX(\ks\partial_y\varphi,\bq)&-&\As u
&=&\mathcal{R}(u,\ks\partial_y\phi,\bq),
\end{array}
\end{equation}
with $\As$ the linear differential operator studied in previous sections, defined in \eqref{As} and that corresponds to the linearisation of the Saint Venant equations about the steady roll-wave $(\Hs,\Qs)$. Recall that
$$
\displaystyle
\As (\hh,\qq)=\left(\begin{array}{c}\displaystyle\ks\partial_y(\cs \hh-\qq)\\
\displaystyle\partial_h\GGs[\hh]+\partial_q\GGs[\qq]-\oms\partial_y \qq\end{array}\right).
$$
Operator $\BWT_0$ is given by
$$
\displaystyle
\BWT_0(\varphi,\bq)=\left(\begin{array}{c}\displaystyle 
\frac{\varphi}{\ks}\Hs'+\d\Hs(\partial_y\varphi,\bq)\\[1em]
\displaystyle 
\frac{\varphi}{\ks}\Qs'+\d\Qs(\partial_y\varphi,\bq)
\end{array}\right),
$$
and $\BWX$ by the fact that $\BWX(k,\bq)$ is
$$
\left(\!\!\begin{array}{c}\displaystyle \ks\partial_y\left((\cs\d\Hs-\d\Qs)(k,\bq)\right)\\
\displaystyle \!\!
\partial_h\GGs[\d\Hs(k,\bq)]+\partial_q\GGs[\d\Qs(k,\bq)]+\partial_k\GGs k+\delta\ks\Qs'\partial_y k
-\oms\partial_y(\d\Qs(k,\bq))\!\!
\end{array}\!\!\right)
$$
and thus may also be written as
$$
\left(\!\!\begin{array}{c}\displaystyle 
-\frac{1}{\ks}\Hs'\,\ks^2 \d\cs(k,\bq)-\Hs\,\ks \d\cs(\partial_y k,\partial_y\bq)+\ks\partial_y\bq\\
\displaystyle
-\frac{1}{\ks}\Qs'\,\ks^2\d\cs(k,\bq)+\BBWX(\partial_yk,\partial_y\bq)
\end{array}\right)
$$
where $\BBWX$ is some differential operator with $1$-periodic coefficients. At last, $\BWT_1(u,k,\bq)$ is a linear differential operator whose coefficients depend at least linearly on $(u,k,\bq)$ and $\mathcal{R}$ is a nonlinear operator acting at least quadratically, informally
$$
\displaystyle
\BWT_1(u,k,\bq)=\mathcal{O}(|k|+|\bq|+|u|),\quad \mathcal{R}(u,k,\bq)=\mathcal{O}(|k|^2+|\bq|^2+|u|^2).
$$

We will split equation \eqref{stv2w1} by projecting it on low Floquet-number critical modes of $\As$. Therefore we first introduce a projection on critical modes of $\As$ for small enough Floquet numbers. 

First recall that, for $l\in[-\pi, \pi]$, $\Asc(l)$ is the operator defined by
$$
\displaystyle
[\Asc(l)(f)](y)=e^{-il y}[\As (e^{il \cdot}f(\cdot))](y)\ 
$$
so that
$$
\widecheck{\As f}(\cdot,l)\ =\ \Asc(l)\check{f}(\cdot,l)\ .
$$
A positive $l_1\in]0,\pi]$ can be chosen small enough so that for any $l$ such that $|l|\leq l_1$ the spectrum  of $\Asc(l)$ in a small $0$-centered ball is given by two spectral curves $\lambda_j(l)$, $j=1,2$, studied in previous section, and may be defined the associated $\Asc(l)$-invariant spectral projection
\begin{equation}\label{Qcc}
\displaystyle
\Qcc(l)=\frac{1}{2\pi i}\int_{\Gamma}(\lambda-\Asc(l))^{-1}d\lambda\ ,
\end{equation}
$\Gamma$ being the boundary of our neighbourhood of $0$.

We further choose a non-increasing $\mathcal{C}^{\infty}$ cut-off function, $\chi: \R\to [0, 1]$, so that, for $l\in\R$,
\begin{equation}\label{cut-off}
\displaystyle
\chi(l)=\left\{\begin{array}{ll}\displaystyle 1\quad{\rm if}\quad|l|\leq 1\\
                                                \displaystyle 0\quad{\rm if}\quad|l|\geq 2,
\end{array}\right.
\end{equation}
and define (with a slight abuse of notation), for $l\in[-\pi, \pi]$, truncated projections
\begin{equation}\label{FSMF}
\begin{array}{ll}
\displaystyle
\Pcc(l)=\Qcc(l)\ \chi\left(\frac{2l}{l_1}\right),\quad
\Psc(l)=\un-\Qcc(l)\ \chi\left(\frac{16l}{l_1}\right),\\
\displaystyle
\PcFSc(l)=\Qcc(l)\ \chi\left(\frac{4l}{l_1}\right),\quad
\PsFSc(l)=\un-\PcFSc(l)\\
\displaystyle
\PcMFc(l)=\Qcc(l)\ \chi\left(\frac{8l}{l_1}\right),\quad
\PsMFc(l)=\un-\PcMFc(l).
\end{array}
\end{equation}
To these operator-valued Bloch-symbols, we associate operators (and use an obvious notation for them). These operators commute and
\begin{equation}\label{proj}
\displaystyle
(\un-\Ps)\PsMF\ =\ 0,\quad (\un-\Ps)\PsFS\ =\ 0.
\end{equation}

We now replace system \eqref{stv2w1} with
{\setlength\arraycolsep{1pt}
\begin{equation}\label{sep}
\left\{\begin{array}{rcl}
\displaystyle
\Big[\PcFS\BWT_0&+&\PcMF\BWT_1(u,\ks\partial_y\varphi,\bq)\Big](\ks\partial_t \varphi,\partial_t \bq)\\[1ex]
&=&\PcFS\BWX(\ks\partial_y\varphi,\bq)+\PcMF\mathcal{R}(u,\ks\partial_y\varphi,\bq)\\[1ex]
\displaystyle
\partial_t u-\As u
&=&\PsFS\BWX(\ks\partial_y\varphi,\bq)+\PsMF\mathcal{R}(u,\ks\partial_y\varphi,\bq)\\[1ex]
\displaystyle
&-&\left[\PsFS\BWT_0+\PsMF\BWT_1(u,\ks\partial_y\varphi,\bq)\right](\ks\partial_t \varphi,\partial_t \bq)
\end{array}\right.
\end{equation}}
supplemented with constraints
\begin{equation}\label{condPS}
\displaystyle
(\un-\Ps)\ u\ =\ 0
\end{equation}
and 
\begin{equation}\label{supp_phi}
\displaystyle
{\rm supp}\big(\mathcal{F}(\varphi),\mathcal{F}(\bq)\big)\ 
\subset\ \left\{l\ \middle|\ \chi\left(\frac{4l}{l_1}\right)=1\right\}.
\end{equation}
Obviously any solution to our new formulation of the problem does provide us with a solution to \eqref{stv2w1}.

A straightforward consequence of \eqref{proj} is that we only need to check that constraint \eqref{condPS} is satisfied at $t=0$. In order to check whether it is so also for assumption \eqref{supp_phi}, we need to describe precisely $\PcFS \BWT_0$. Of course we also need this detailed description to compare first equation of system \eqref{sep} with the Whitham's system.

But our spectral study of previous section (see Lemmas~\ref{eigenvect0} and~\ref{eigenvect1}) does give us the needed expansion of $\Qcc$. Correspondingly, we can write eigenprojection $\PcFS$ as
\begin{equation*}
\displaystyle
\PcFSc\check{u}(y,l)=\chi\left(\frac{4l}{l_1}\right)\left(<\tilde{v}_1(\cdot,l),\check{u}(\cdot,l)>v_1(y,l)+<\tilde{v}_2(\cdot,l),\check{u}(\cdot,l)>v_2(y,l)\right)
\end{equation*}
where again $<\cdot,\cdot>$ is the scalar product of $L^2_{\txt{per}}$.

Besides, one readily obtains that, for any pair $(\varphi,\bq)$ such that \eqref{supp_phi} is satisfied,
\begin{equation}
\mathcal{J}[\BWT_0(\ks\varphi,\bq)](l)=\ks\hat{\varphi}(l)(v_1(l)+\mathcal{O}(l^2))+\hat{\bq}(l)(\partial_{\bq}\Ms v_2(l)+\mathcal{O}(l))
\end{equation}
therefore (using also $0=<\partial_l\tilde{v}_2(0),v_1(0)>+<\tilde{v}_2(0),\partial_l v_1(0)>$)
\begin{equation}\label{dt_phi_q}
\begin{array}{ll}
\displaystyle
<\tilde{v}_1(l),\mathcal{J}[\BWT_0(\ks\varphi,\bq)](l)>=&
\left(1+\widehat{\widetilde{\mathcal{B}}}^T_{k,2}(l)\right)\ks\hat{\phi}(l)+\widehat{\widetilde{\mathcal{B}}}^T_{k,1}(l)\ \hat{\bq}(l),\\
\displaystyle
<\tilde{v}_2(l),\mathcal{J}[\BWT_0(\ks\varphi,\bq)](l)>=&\left(\d\Ms+\widehat{\widetilde{\mathcal{B}}}^T_{\bq,1}(l)\right)[\widehat{\ks\partial_y\varphi}(l),\hat{\bq}(l)]\\
&-i<\partial_l\tilde{v}_2(0),v_1(0)>\widehat{\ks\partial_y\varphi}(l),
\end{array}
\end{equation}
with formally $\|\widehat{\widetilde{\mathcal{B}}}^T_{k,1}(l)\|,\|\widehat{\widetilde{\mathcal{B}}}^T_{\bq,1}(l)\|=\mathcal{O}(l)$ and $\|\widehat{\widetilde{\mathcal{B}}}^T_{k,2}(l)\|=\mathcal{O}(l^2)$.

Likewise, \eqref{supp_phi} implies
\begin{equation*}\label{dx_phi_q}
\begin{array}{ll}
\displaystyle
<\tilde{v}_1(l),\mathcal{J}[\BWX(k,\bq)](l)>=&-\ks^2\d\cs(\hat{k}(l),\hat{\bq}(l))+\widehat{\widetilde{\mathcal{B}}}^X_{k,1}(l)[\hat{k}(l),\hat{\bq}(l)],\\
\displaystyle
<\tilde{v}_2(l),\mathcal{J}[\BWX(k,\bq)](l)>=&-\Ms\ks\d\cs(\widehat{\partial_y k}(l),\widehat{\partial_y\bq}(l))+\ks\widehat{\partial_y\bq}(l)\\
\displaystyle
&+i<\partial_l\tilde{v}_2(0),v_1(0)>\ks^2\d\cs(\widehat{\partial_y k}(l),\widehat{\partial_y\bq}(l))\\
\displaystyle
&+\widehat{\widetilde{\mathcal{B}}}^X_{\bq,2}(l)[\hat{k}(l),\hat{\bq}(l)],
\end{array}
\end{equation*}
with  $\|\widehat{\widetilde{\mathcal{B}}}^X_{k,1}(l)\|=\mathcal{O}(l)$ and $\|\widehat{\widetilde{\mathcal{B}}}^X_{\bq,2}(l)\|=\mathcal{O}(l^2)$.

Now, from definition of $\PcMF$ is derived for any function f
$$
\displaystyle
{\rm supp}_l\mathcal{J}[\PcMF f](\cdot,l)\subset\left\{l\,\middle|\,\chi\left(\frac{4l}{l_1}\right)=1\right\}.
$$
Since $\left[\begin{array}{cc}1&0\\\partial_k\Ms&\partial_{\bq}\Ms\end{array}\right]$ is invertible, as a result, assumption \eqref{supp_phi} is seen to be satisfied whenever it is at $t=0$. From now on we do not repeat but always assume \eqref{supp_phi} is satisfied.

We still need to relate the first equation of system \eqref{sep} with the Whitham's system but it is now a straightforward task. Let us split $\BWT_1(u,k,\bq)$ into 
$$
\BWT_1(u,k,\bq)(\varphi,\tilde{\bq})=\BWT_{1,\varphi}(u,k,\bq)\varphi+\BWT_{1,\bq}(u,k,\bq)\tilde{\bq}
$$
and define $\pi_j$ by $\widehat{\pi_j f}(l)=<\tilde{v}_j(\cdot,l),\check{f}(\cdot,l)>$. Setting $k=\ks\partial_y\varphi$ and applying $\pi_1$ to the first line of \eqref{sep} leads to 
\begin{equation*}
\displaystyle
\label{eq_phi}
\begin{array}{rcl}
\Big(1+\BBWT_{k,2}&+&\pi_1\PcMF\BWT_{1,\varphi}(u,k,\bq)\Big)\partial_t\left(\ks\varphi\right)\\
&=&-\ks^2\d\cs(k,\bq)+\BBWX_{k,1}(k,\bq)+\pi_1\PcMF\mathcal{R}(u,k,\bq)\\
&&-\pi_1\PcMF\BWT_{1,\bq}(u,k,\bq)\partial_t\bq-\BBWT_{k,1}\partial_t\bq\ .
\end{array}
\end{equation*}
If $l_1$ is small enough, then, as long as $(k,\bq,u)$ is kept small enough, this can be turned into
\begin{equation*}
\displaystyle
\begin{array}{rcl}
\partial_t \left(\ks\varphi\right)&+&\ks^2\d\cs(k,\bq)\\
&=&-\left[\left(1+\BBWT_{k,2}+\pi_1\PcMF\BWT_{1,\varphi}(u,k,\bq)\right)^{-1}-1\right]\left(\ks^2\d\cs(k,\bq)\right)\\
&+&\left(1+\BBWT_{k,2}+\pi_1\PcMF\BWT_{1,\varphi}(u,k,\bq)\right)^{-1}\left(\BBWX_{k,1}(k,\bq)+\pi_1\PcMF\mathcal{R}(u,k,\bq)\right)\\
&-&\left(1+\BBWT_{k,2}+\pi_1\PcMF\BWT_{1,\varphi}(u,k,\bq)\right)^{-1}\left(\left(\pi_1\PcMF\BWT_{1,\bq}(u,k,\bq)+\BBWT_{k,1}\right)\partial_t\bq\right)
\end{array}
\end{equation*}
and, denoting the right-hand side of the former equation by 
$$
F_k(u,k,\bq)\ -\ \BBWT_{k,\bq}(u,k,\bq)\partial_t\bq\ -\ \BBWX_k(k,\bq),
$$
with formally $\|\widehat{\widetilde{\mathcal{B}}}^T_{k,\bq}(u,k,\bq)(l)\|=\mathcal{O}(|l|+|u|+|k|+|\bq|)$, $\|\widehat{\widetilde{\mathcal{B}}}^X_k\|=\mathcal{O}(l)$ and $\|F_k(u,k,\bq)\|=\mathcal{O}(|u|^2+|k|^2+|\bq|^2)$,
gives
\begin{equation}
\label{eq_prephi1}
\begin{array}{rcl}
\partial_t \left(\ks\varphi\right)&+&\BBWT_{k,\bq}(u,k,\bq)\partial_t\bq\\ 
&+&\ks^2\d\cs(k,\bq)\ +\ \BBWX_k(k,\bq)
\ =\ F_k(u,k,\bq)
\end{array}
\end{equation}
therefore
\begin{equation}
\label{eq_phi1}
\begin{array}{rcl}
\partial_tk&+&\partial_y\left(\BBWT_{k,\bq}(u,k,\bq)\partial_t\bq\right)\\
&+&\ks^2\d\cs(\partial_yk,\partial_y\bq)\,+\,\partial_y\left(\BBWX_k(k,\bq)\right)
\,=\,\partial_y\left(F_k(u,k,\bq)\right).
\end{array}
\end{equation}

Applying $\pi_2$ rather than $\pi_1$ leads to
\begin{equation*}
\displaystyle
\begin{array}{rl}
\left(\d\Ms+\BBWT_{\bq,1}\right)&\left(\partial_t k,\partial_t\bq\right)
+\pi_2\PcMF\BWT_{1,\varphi}(u,k,\bq)\partial_t\left(\ks\varphi\right)\\[1ex]
\displaystyle
&=-\ks\Ms\d\cs(\partial_y k,\partial_y\bq)+\ks\partial_y\bq+\BBWX_{\bq,2}(k,\bq)\\[1ex]
\displaystyle
&+\ \pi_2\PcMF\mathcal{R}(u,k,\bq)-\pi_2\PcMF\BWT_{1,\bq}(u,k,\bq)\partial_t\bq\\[1ex]
\displaystyle
&+\ i<\partial_l\tilde{v}_2(0),v_1(0)>\left[\partial_tk+\ks^2\d\cs(\partial_yk,\partial_y\bq)\right].
\end{array}
\end{equation*}
Using \eqref{eq_prephi1}, the equation may be written as
\begin{equation}
\label{eq_q1}
\begin{array}{rcl}
\Big(\d\Ms&+&\BBWT_{\bq}(u,k,\bq)\Big)\left(\partial_t k,\partial_t\bq\right)\\
&+&\ks\Ms\d\cs(\partial_y k,\partial_y\bq)-\ks\partial_y\bq
+\BBWX_{\bq}(k,\bq)\ =\ F_{\bq}(u,k,\bq)
\end{array}
\end{equation}
with formally $\|\widehat{\widetilde{\mathcal{B}}}^T_{\bq}(u,k,\bq)(l)\|=\mathcal{O}(|l|+|u|+|k|+|\bq|)$, $\|\widehat{\widetilde{\mathcal{B}}}^X_{\bq}\|=\mathcal{O}(l^2)$ and $\|F_{\bq}(u,k,\bq)\|=\mathcal{O}(|u|^2+|k|^2+|\bq|^2)$.

Again, under smallness assumptions, one can invert
$$
\displaystyle
\left[\begin{array}{cc}
\displaystyle
1&il\ \widehat{\widetilde{\mathcal{B}}}^T_{k,\bq}(u,k,\bq)(l)\\
\displaystyle
\partial_k\Ms+\widehat{\widetilde{\mathcal{B}}}^T_{\bq,k}(u,k,\bq)(l)
&\partial_{\bq}\Ms+\widehat{\widetilde{\mathcal{B}}}^T_{\bq,\bq}(u,k,\bq)(l)
\end{array}\right].
$$
Thus equations (\ref{eq_phi1},\ref{eq_q1}) yield an evolution system for $(k,\bq)$. As a result, setting $\V= \left((k,\bq), u\right)$, system \eqref{sep} may be written in short as
\begin{equation}\label{sep2}
\displaystyle
\partial_t\V-\Lambda\V=\NL(\V),
\end{equation}
where $\Lambda$ is a lower triangular linear operator written in Fourier-Bloch variables (Fourier in $(k,\bq)$, Bloch in $u$) as\footnote{With a slight abuse of notations in the use of $\ \widecheck{\ }$.}
$$
\Lc(l)=
\left(\begin{array}{cc} \displaystyle \Acf(l) & 0\\
\displaystyle \fc{b}\,(l) &\Asc(l)
\end{array}\right), 
$$
with $\fc{b}(l)=\mathcal{O}(l)$ a bounded operator (from $\C^2$ to $1$-periodic functions) and $\Acf(l)\in \mathcal{M}_2(\C)$ such that $\Acf(l)=\mathcal{O}(l)$, and $\NL$ a nonlinear operator such that $\NL(\V)=\mathcal{O}(|\V|^2)$. From the derivation of the critical system, the spectrum of $\Acf(l)$ is easily seen to be given by the values $\lambda_j(il)$, $j=1,2$, of the previously introduced spectral curves $\lambda_j$, $j=1,2$.

In order to emphasize the mode separation, we write $\V=(\vc,\vs)$ with $\vc=(k,\bq)$ and $\vs=u=(\hh,\qq)$. Then system \eqref{sep2} may be written
\begin{equation}\label{diag}
\left\{
\begin{array}{rcl}
\displaystyle
\partial_t \vc\,-\,\Ac\,\vc&=&\NLc(\vc,\vs)\\
\displaystyle
\partial_t \vs\,-\,\As\,\vs-b \vc&=&\NLs(\vc,\vs)
\end{array}
\right.
\end{equation}
where $\NLc(\vc,\vs)$, $\NLs(\vc,\vs)=\mathcal{O}(|\vc|^2+|\vs|^2)$ and $\NLc$ is of the following form
$$
\NLc(\vc,\vs)\ =\ f(\vc,\vs)\,\partial_y F(\vc,\vs)\ +\ \rho\,\NNLc(\vc,\vs)
$$
with $f(\vc,\vs)$, $F(\vc,\vs)=\mathcal{O}(|\vc|+|\vs|)$, $\NNLc(\vc,\vs)=\mathcal{O}(|\vc|^2+|\vs|^2)$ and  $\rho$ is such that $\widehat{\rho}(l)=\mathcal{O}(l)$.  This former fact is trivial for the part coming from $\pi_1$ through equation \eqref{eq_phi1}. As for the contribution of $\pi_2$, it follows from the following fact:
$$
\widehat{\pi_2 (h,q)}(l)\ =\ <\widecheck{h}(\cdot,l)>\ +\ \mathcal{O}(l)[(h,q)]
\ =\ \widehat{h}(l)\ +\ \mathcal{O}(l)[(h,q)]\ .
$$

This is the end of our preparation of the system and we can now build a family of approximate solutions according to the desired \emph{ansatz} and then achieve the proof with the construction of a family of solutions close to our family of approximate solutions. For this purpose we are now in such a position that we can follow the strategy explained in \cite{D3S} and therefore we will mostly sketch the end of the proof.

\subsection{Approximate solutions}

We now fix $a>0$, $m\geq3$, and $T_0>0$ small enough and consider a smooth solution $(k,\bq)$ to system \eqref{just_whit} in $L^\infty\left([0,T_0];\mathcal{Y}^a_0\right)$. For further implicit uses, note that for some $C(a,m)$ stands
$$
\xnorm{\ \cdot\ }{a}{m}\ \leq\ C(a,m)\,\ynorm{\ \cdot\ }{a}{0}
$$
and that, for any $\varepsilon>0$, $\ynorm{\,\cdot\,}{a}{0}$ is turned into $\ynorm{\,\cdot\,}{\varepsilon a}{0}$ by the transformation $f\mapsto f(\varepsilon\,\cdot)$ so that in particular $\ynorm{\,\cdot\,}{a}{0}$ is invariant under this transformation whenever $\varepsilon\leq1$.

In order to build approximate solutions in the long-wavelength regime we associate to long-wavelength profiles $(\Vc,\Vs)$ and any $\varepsilon>0$ $\varepsilon$-residuals $\Resc(\Vc,\Vs)$ and $\Ress(\Vc,\Vs)$ through
$$
\begin{array}{rcl}
\displaystyle
\Resc(\Vc,\Vs)(X,T)&=&\left[\partial_t \vc-\Ac \vc-\NLc(\vc,\vs)\right]\left(\frac{X}{\varepsilon},\frac{T}{\varepsilon}\right)\\
\displaystyle
\Ress(\Vc,\Vs)(X,T)&=&\left[\partial_t \vs-\As \vs-b\vc-\NLs(\vc,\vs)\right]\left(\frac{X}{\varepsilon},\frac{T}{\varepsilon}\right)
\end{array}
$$
where $(\vc,\vs)$ is defined by
$$
\vc(y,t)\ =\ \Vc(\varepsilon y,\varepsilon t)\ ,\quad
\vs(y,t)\ =\ \Vs(\varepsilon y,\varepsilon t)\ .
$$
Obviously, for a given $\varepsilon>0$, the above $(v^c,v^s)$ is a solution to \eqref{diag} if and only if $\Resc(\Vc,\Vs)$ and $\Ress(\Vc,\Vs)$ vanish.

The next proposition provides us with the needed approximate solutions. Yet, to be able to prove it we need to understand the behaviour of operators, defined in Bloch variables, with respect to dilatation. For this purpose, let us denote $D_\alpha$ the dilatation operator, $D_{\alpha} (f)=f(\alpha\,\cdot)$.

First recall the diagonalisation formula
$$
[Tf](x)\ =\  \frac{1}{\sqrt{2\pi}}\int_{-\pi}^\pi e^{ilx}[\widecheck{T}(l)\check{f}(\cdot,l)](x)\d l
$$
where
$$
[\widecheck{T}(l)g](x)\ =\ e^{-ilx}[T(e^{il\cdot}g(\cdot))](x)\ .
$$
From this one deduces 
$$
[D_{\varepsilon^{-1}}TD_{\varepsilon}](f)(x)
\ =\ 
\frac{1}{\sqrt{2\pi}}\int_{-\pi}^\pi e^{ilx}[D_{\varepsilon^{-1}}\widecheck{T}(\varepsilon l)D_{\varepsilon}](\check{f}(\cdot,l))(x)\d l\ .
$$
Note that one can not infer form this formula a Bloch transform since the periodicity may be lost in the process. For instance, in the case where $\widecheck{T}(l)=P(l,y,\partial_y)$ with $P$ a symbol $1$-periodic in $y$, in the formula appears $P\left(\varepsilon l,\frac{y}{\varepsilon},\varepsilon \partial_y\right)$. However note that even in this case thanks to the original periodicity the $y/\varepsilon$-dependency is almost harmless.

As for Bloch-Fourier operators, let us look at an operator defined through
$$
\widehat{Tf}(l)\ =\ \int_0^1 [\cf{\tau}(l)\widecheck{f}(\cdot,l)](y)\d y
$$
(with an extended definition of the Bloch transform). Then formally
$$
\mathcal{F}\left([D_{\varepsilon^{-1}}TD_{\varepsilon}](f)\right)(l)
\ =\ \sum_{j\in\Z}\widehat{f}\left(l+\frac{2\pi j}{\varepsilon}\right) \int_0^1 [\cf{\tau}(\varepsilon l)[e^{2i\pi j\cdot}]](y)\d y
$$
which, when $\cf{\tau}(l)=P(l,y,\partial_y)$, turns into
$$
\sum_{j\in\Z}\widehat{f}\left(l+\frac{2\pi j}{\varepsilon}\right) \int_0^{\varepsilon} \left[P\left(\varepsilon l,\frac{y}{\varepsilon},\varepsilon \partial_y\right)\left[e^{2i\pi \frac{j}{\varepsilon}\cdot}\right]\right](y)\frac{\d y}{\varepsilon}\ .
$$
Note that, when $f$ is low-frequency, it only involves bounded $j/\varepsilon$ so that, when moreover $P$ is $1$-periodic in $y$, again oscillations are harmless.

At last, let us consider an operator defined through
$$
[Tf](x)\ =\ \frac{1}{\sqrt{2\pi}}\int_{-\pi}^\pi e^{ilx}[\fc{\tau}(x)\hat{f}(\cdot)](l)\d l
$$
with $\fc{\tau}(\cdot)$ $1$-periodic. Then
$$
\left[D_{\varepsilon^{-1}}TD_{\varepsilon}\right](f)(x)
\ =\ \frac{1}{\sqrt{2\pi}}\int_{-\frac{\pi}{\varepsilon}}^{\frac{\pi}{\varepsilon}} e^{ilx}\left[\left[D_{\varepsilon}\fc{\tau}\left(\frac{x}{\varepsilon}\right)D_{\varepsilon^{-1}}\right](\widehat{f}\,)\right](l)\d l
$$
which, when $\fc{\tau}(y)=P(l,y)$, turns into
$$
\frac{1}{\sqrt{2\pi}}\int_{-\frac{\pi}{\varepsilon}}^{\frac{\pi}{\varepsilon}} e^{ilx}\,P\left(\varepsilon l,\frac{x}{\varepsilon}\right)\,\widehat{f}(l)\d l\ .
$$
Again, when $f$ is low-frequency, oscillations (here in $\pi/\varepsilon$ and $x/\varepsilon$) are harmless.

We are now in position to state the following proposition.

\begin{proposition}\label{approximate}
Let $m\geq 3$, $M\geq1$ and $a_0>0$. There exists $\eta_1>0$, $\varepsilon_1>0$ and a constant $K_1$ such that, for any $0<a\leq a_0$, if $(k,\bq)$ is a solution to the Whitham system such that $\sup_{T\in[0,T_0]}\ynorm{(k,\bq)(\cdot,T)}{a}{m}\leq \eta_1$, then, for any $0<\varepsilon<\varepsilon_1$, there exists $(\Vceps,\Vseps)$ such that 
{\setlength\arraycolsep{1pt}
\begin{eqnarray}
\displaystyle
\sup_{T\in[0, T_0]}\xnorm{\Vceps(\varepsilon\cdot,T)-(k,\bq)(\varepsilon \cdot,T)}{a/\varepsilon}{m}
&\leq& K_1\left[\varepsilon+\sup_{T\in[0, T_0]}\hulnorm{(k,\bq)(\cdot,T)-(\ks,\bqs)}{m}^2\right],\nonumber\\
\displaystyle
\sup_{T\in[0, T_0]}\xnorm{\Vceps(\varepsilon\cdot,T)}{a/\varepsilon}{m}&\leq& K_1,\nonumber\\
\displaystyle
\sup_{T\in[0, T_0]}\xnorm{\Vseps(\varepsilon\cdot,T)}{a/\varepsilon}{m}&\leq& K_1\sup_{T\in[0, T_0]}\hulnorm{(k,\bq)(\cdot,T)-(\ks,\bqs)}{m}^2,\nonumber\\
\displaystyle
\sup_{T\in[0, T_0]}\xnorm{\Resc(\Vceps,\Vseps)(\varepsilon\cdot,T)}{a/\varepsilon}{m}&\leq&K_1\varepsilon^M,\nonumber\\
\displaystyle
\sup_{T\in[0, T_0]}\xnorm{\Resc(\Vceps,\Vseps)(\varepsilon\cdot,T)}{a/\varepsilon}{m}&\leq&K_1\varepsilon^M.\nonumber
\end{eqnarray}}
\end{proposition} 

To prove Proposition~\ref{approximate}, one search for $(\Vceps,\Vseps)$ in the form
$$
\begin{array}{rcl}
\displaystyle
\Vceps&=&\Vc_0\ +\ \varepsilon\,\Vc_1\ +\ \cdots\ \varepsilon^M\,\Vc_{M}\\
\displaystyle
\Vseps&=&\Vs_0\ +\ \varepsilon\,\Vs_1\ +\ \cdots\ \varepsilon^M\,\Vs_{M}\\
\end{array}
$$
and obtain a hierarchy of equations for $(\Vc_j,\Vs_j)_{0\leq j\leq M}$.

At step $j$, $\Vs_j$ is obtained as a function of $(\Vc_i,\Vs_i)_{0\leq i<j}$ and $\Vc_j$ by solving an equation (not of evolution type) through the implicit function theorem thanks to the inversibility of $\As\PsFS$. For $j\neq0$ this function is linear in $\Vc_j$. Using this expression for $\Vs_j$ one obtains an evolution-type equation for $\Vc_j$ in terms of $(\Vc_i,\Vs_i)_{0\leq i<j}$. For $j=0$, the equation for $\Vc_0$ shares its linearization around $(\ks,\bqs)$ with the equation for $(k,\bq)$. For $j\neq0$, the equation is \emph{linear} hyperbolic.

Note that the intricated form of the proposition is a consequence of the fact that the norm $\xnorm{\cdot}{a}{m}$ badly scales. Note also that, although, up to the expression of $\Vs_j$, the proof of the proposition follows the lines of the formal derivation of the Whitham system, the two expansions may differ even though they share the same starting point $(k,\bq)$ at the linear level.

\subsection{From residuals to reminders}

We now look for a family $((\vceps,\vseps))_{0<\varepsilon<\varepsilon_1}$ of solutions to system \eqref{diag} in the form
$$
\displaystyle 
(\vceps,\vseps)(x,t)\ =\ (\Vceps,\Vseps)(\varepsilon x,\varepsilon t)\ +\ \varepsilon^M (\rceps,\rseps)(x,t)
$$
with $(\rceps,\rseps)$ uniformly bounded on $[0, T_1/\varepsilon]$ (where $T_1$ is some fixed time $0<T_1<T_0$).

Substituting this ansatz into \eqref{diag} yields an equation we write as
\begin{equation}\label{diag_rest}
\left\{\begin{array}{ll}
\displaystyle
\partial_t \rceps\ -\ \Ac \rceps\ =\ \Nceps(\rceps,\rseps)\\
\displaystyle
\partial_t \rseps\ -\ \As \rseps\ =\ \Nseps(\rceps,\rseps)
\end{array}\right.\ .
\end{equation}
Note that a linear term $b\,\rceps$ has been put in the right-hand side of the second equation of system \eqref{diag_rest} so as to deal with a diagonal form in the left-hand part. As usual this family of systems are solved by a fix point argument.

For this purpose, we report the following estimates
$$
\begin{array}{rcl}
\displaystyle
\xnorm{\Nceps(\rceps,\rseps)}{a}{m}
&\leq&
K_1\ +\ K_{\eta_1,l_1}\left(\xnorm{\rceps}{a}{m}+\xnorm{\rseps}{a}{m}\right)
\ +\ \varepsilon^M K(R_{\txt{c}},R_{\txt{s}})\\
\displaystyle
\xnorm{\Nseps(\rceps,\rseps)}{a}{m-2}
&\leq&
K_1\ +\ K_{\eta_1,l_1}\left(\xnorm{\rceps}{a}{m}+\xnorm{\rseps}{a}{m}\right)
\ +\ \varepsilon^M K(R_{\txt{c}},R_{\txt{s}}),
\end{array}
$$
valid for all $0<a<\frac{a_0}{\varepsilon}$ whenever $\xnorm{\rceps}{a}{m}\leq R_{\txt{c}}$ and $\xnorm{\rseps}{a}{m}\leq R_{\txt{s}}$, where $K_{\eta_1,l_1}$ goes to zero when $(\eta_1,l_1)$ goes to zero. Recall that $l_1$ is a cut-off parameter and $\eta_1$ is a size parameter for $(k,\bq)$. Actually $K_1$ also (badly) depends on $l_1$ but this latter point is not prejudicial.

It is crucial to note that through these estimates one undergoes a loss of derivatives. This is a consequence of the fact that our initial nonlinear change of variables has turned the semilinear Saint-Venant equations into quasilinear ones. This former point also explains why it would not be harder to deal with a more physical viscosity. Unfortunately, at this stage, one can not exploit any smoothing coming from the linear operator $\As\PsFS$. But this is a common fact that in the process of justifying an ansatz one usually loses something. Here instead of losing regularity we will choose to lose analyticity. We will establish estimates with a width of analyticity descreasing at a linear pace, so that we will work on a finite time, even shorter than $T_0$.

For this we introduce a \emph{smoothing} operator as follows. Let us fix $K_0>0$. We use the fact that the spectrum of $\As$ is of \emph{upper bounded} real part to choose $K'_0$ sufficiently large so that, for any $l\in[-\pi,\pi]$, the real part of the spectrum of $\Asc(l)-K'_0|l|$ is upper bounded by $-K_0|l|$. We further define operators $k'_0$ and $\mathcal{S}_\varepsilon(t)$ through their Bloch symbols $\widecheck{k'_0}(l)=K'_0|l|$ and $\widecheck{\mathcal{S}_\varepsilon(t)}(l)=e^{\left(\frac{a_0}{\varepsilon}-K'_0 t\right)|l|}$. 

Then we set $\displaystyle \left(\widetilde{\rceps},\widetilde{\rseps}\right)(t)=\left(\mathcal{S}_\varepsilon(t)\rceps(t),\mathcal{S}_\varepsilon(t)\rceps(t)\right)$ for times $t$ satisfying
$$
0\ \leq\ t\ <\ \frac{a_0}{K'_0\,\varepsilon}\ .
$$ 
Estimates are thus established now in $\xnorm{\cdot}{0}{m}$ norms. The system for the evolution of $\left(\widetilde{\rceps},\widetilde{\rseps}\right)$ is similar to the previous one, with the same kind of estimates for nonlinear terms, but with linear operator $(\Ac,\As)$ replaced with $(\Ac-k'_0,\As-k'_0)$. The reason for this change of unknowns is that now one can prove for $(\Ac-k'_0,\As-k'_0)$ a maximal regularity result in $\mathcal{X}_m^0$-valued H\"older spaces, readily similar to Lemma 6.3 in \cite{D3S}. This was the missing part to close in a classical way a fix-point iteration scheme.

The last thing we should say is that $\varphi_{0,\varepsilon}$ is then recovered by integrating over time $T_1/\varepsilon$ equation \eqref{eq_prephi1}.


\section{Conclusion}

In this paper, we derived formally first-order and second-order averaged equations for shallow water flows that describes the dynamics of modulated roll-waves and provided two set of justification results. 

On the one hand, we carried out a spectral stability analysis of roll-waves using Bloch transform in the regime of small wavenumber perturbations. We first related the index of stability of viscous roll-waves with the hyperbolicity of the first-order Whitham's equations just as it was done by Serre for general viscous conservation laws \cite{Serre} or Johnson, Zumbrun and  Bronski for generalized Korteweg-de Vries equations \cite{Z_kdv}. However, in both latter papers, only the Evans function framework was used. Here we proved the Bloch transform framework is a more natural tool, by extending such stability analysises into two directions : we relate not only eigenvalues but also eigenvectors and we relate the stability of steady solutions to the \emph{second-order} Whitham's equations with the \emph{parabolicity} of spectral curves at the origin. 

On the other hand, in the spirit of what has been done for reaction-diffusion equations by Doelman, Sandstede, Scheel and Schneider \cite{D3S}, we justified rigorously the inviscid Whitham's equations in the natural hyperbolic scaling. More precisely, we proved that, given a solution to the inviscid Whitham's system, there exist solutions to the viscous shallow water equations on asymptotically large time that are close to modulated roll-waves whose first-order expansion is described by our solution to the Whitham's system at the linear level. This justification is performed under weak stability assumptions : at the origin, tangency to the imaginary axis of spectral curves. From numerical investigations, it seems that it will enable us to apply this nonlinear justification up to the limiting homoclinic travelling waves (whose spectrum yields unstablility but is tangent to to the imaginary axis \cite{Ba_Jo_Ro_Zu}). This weak stability assumption has a counterpart in the required analyticity of solutions.

But probably the main flaw of this justification is that our solution to the Whitham's system describes the first order of the roll-wave profile only \emph{at the linear level}. This is a consequence of both the $\Yphi$ change of variables and the hyperbolic scaling. Such an issue would not occur with a diffusive scaling. Yet the counterpart would be that the size of the allowed perturbations of $(k,\bq)$ would not be anymore $\mathcal{O}(1)$ but $\mathcal{O}(\varepsilon)$ (recall $\varepsilon>0$ is the characteristic wavenumber of the modulation).

Moreover the nonlinear justifications need some regularity so that, with the hyperbolic scaling, there is little hope to construct modulated roll-waves converging at $\pm\infty$ to roll-waves with different wavenumbers and local discharge rate since they would correspond to shocks in parameters $(k,\bq)$. However with a parabolic scaling on may expect to justify viscous shocks of the second-order modulation system as modulated roll-wave profile for solutions to the Saint-Venant system.

At last a second-order justification would probably also enlarge the time validity of the approximation from $[0,T_0/\varepsilon]$ with $T_0$ some time imposed by the equations to $[0,T/\varepsilon^2]$ for any fixed time $T$. 

For all these reasons, a natural direction would be now to justify a second-order modulation in a parabolic scaling as it was also done in \cite{D3S} for reaction-diffusion equations. 

\appendix

\section{Whitham's dispersion in a geometric way}

In this appendix, we explain how to fill the small gap between analysis in subsection \ref{subsec_Evans} and the one in \cite{Serre} (for different equations). We just need to replace our straigthforward computations with more geometric ones. Yet we keep the framework of subsection \ref{subsec_Evans} : frame $(x-\cs t,t)$, $(c,\bq)$-parametrization...

We turn $2\times2$-determinant giving Whitham's dispersion into a $3\times3$-determinant. This is done interpreting the set of periodic travelling-wave solutions as a submanifold.

The tangent space of the set of periodic travelling-wave solutions (indentified when being equals up to translation) at $\Hs$ is
$$
\left\{ \beta_1\partial_c\Hs\ +\ \beta_2 h_2^0\ +\ \beta_3\partial_{\bq}\Hs\ \middle|\ 
(\beta,\gamma)\in \ker Z
\right\}
$$
where the existence of such a $h_2^0$ is provided by the analysis of \eqref{profil} and 
$Z$ is a linear operator from $\R^3\times\R$ to $\R^2$ defined by
\begin{eqnarray*}
Z_1(\beta,\gamma)&=& 
\beta_1[\partial_c\Hs]+\beta_2 [h_2^0]+\beta_3[\partial_{\bq}\Hs]+\gamma\Hs'(0)\ ,\\
Z_2(\beta,\gamma)&=& 
\beta_1[\partial_c\Hs']+\beta_2 [h_2^0{}']+\beta_3[\partial_{\bq}\Hs']+\gamma\Hs''(0)\ .
\end{eqnarray*}

Note that on $\ker Z$
\begin{equation*}
\begin{array}{ll}
\d \Ls(\beta,\gamma)\ =\ \gamma\ ,&\d \cs(\beta,\gamma)\ =\ \beta_1\ ,\\
\d \bqs(\beta,\gamma)\ =\ \beta_3\ ,&
\d\Hs(\beta,\gamma)\ =\ \beta_1\partial_c\Hs\ +\ \beta_2 h_2^0\ +\ \beta_3\partial_{\bq}\Hs\ ,
\end{array}
\end{equation*}
and
\begin{equation*}
\d\Ms(\beta,\gamma)\,=\,\gamma\ks[\Hs(0)-\Ms]+<\d\Hs(\beta,\gamma)>_{\Ls}.
\end{equation*}
Fix $(\lambda,\nu)\in\C^2$ and define $T(\lambda,\nu)$ a linear operator from $\R^3\times\R$ to $\R^2$ by
\begin{equation*}
\begin{array}{l}
T_1(\beta,\gamma)=
-\lambda\ks^2\gamma+\nu\ks^2\beta_1\ ,\\
T_2(\beta,\gamma)= 
\lambda\left[\beta_1<\partial_c\Hs>_{\Ls}+\beta_2<h_2^0>_{\Ls}+\beta_3<\partial_{\bq}\Hs>_{\Ls}\right]+\nu\ks\left[\Ms\beta_1-\beta_3\right].
\end{array}
\end{equation*}
Operator $T(\lambda,\nu)$ co\"incide with $(\lambda\d\ks+\nu\ks\d\cs,\lambda\d\Ms+\nu\ks(\ks\d\cs-\d\bqs))$ on $\ker Z$, since as in subsection \ref{subsec_Evans} we have imposed $\Ms=\Hs(0)$.

Now, up to some non-zero $\Gamma$, $\Gamma'$,
$$
D(\lambda,\nu)\ =\ \Gamma\det (T(\lambda,\nu)_{|\ker Z})\ =\ \Gamma'\det((T(\lambda,\nu),Z))
$$
and $\det((T(\lambda,\nu),Z))$ equals
$$
\left|\begin{array}{cccc}
\nu\ks^2&0&0&-\lambda\ks^2\\
\lambda<\partial_c\Hs>_{\Ls}+\nu\ks\Ms&\lambda<h_2^0>_{\Ls}&\lambda<\partial_{\bq}\Hs>_{\Ls}-\ks\nu&0\\
\left[\partial_c \Hs\right]&[h_2^0]&[\partial_{\bq}\Hs]&\Hs'(0)\\
\left[\partial_c \Hs'\right]&[h_2^0{}']&[\partial_{\bq}\Hs']&\Hs''(0)
\end{array}\right|
$$
and may be reduced to
$$
\ks^2\left|\begin{array}{ccc}
\lambda^2<\partial_c\Hs>_{\Ls}+\lambda\nu\ks\Ms&\lambda<h_2^0>_{\Ls}&<\partial_{\bq}\Hs>_{\Ls}-\nu\\
\lambda\left[\partial_c \Hs\right]+\nu\Hs'(0)&[h_2^0]&[\partial_{\bq}\Hs]\\
\lambda\left[\partial_c \Hs'\right]+\nu\Hs''(0)&[h_2^0{}']&[\partial_{\bq}\Hs']
\end{array}\right|
$$
which easily compares with the main part of the Evans function $E(\lambda,e^\nu)$.


\begin{thebibliography}{10}

\bibitem{Ba_Jo_No_Ro_Zu_2}
Blake {Barker}, Mathew~A. {Johnson}, Pascal {Noble}, L.~Miguel {Rodrigues}, and
  Kevin {Zumbrun}.
\newblock Spectral stability of periodic viscous roll waves.
\newblock In preparation.

\bibitem{Ba_Jo_No_Ro_Zu_1}
Blake {Barker}, Mathew~A. {Johnson}, Pascal {Noble}, L.~Miguel {Rodrigues}, and
  Kevin {Zumbrun}.
\newblock Whitham averaged equations and modulational stability of periodic
  traveling waves of a hyperbolic-parabolic balance law.
\newblock {\em ArXiv e-prints, arXiv:1008.4729v2}, 2010.
\newblock Submitted.

\bibitem{Ba_Jo_Ro_Zu}
Blake {Barker}, Mathew~A. {Johnson}, L.~Miguel {Rodrigues}, and Kevin
  {Zumbrun}.
\newblock {Metastability of solitary roll wave solutions of the St. Venant
  equations with viscosity}.
\newblock {\em ArXiv e-prints, arXiv:1007.5262v1}, 2010.
\newblock Submitted.

\bibitem{Boudlal}
Abdelaziz Boudlal and Valérie~Yu Liapidevskii.
\newblock Stability of regular roll waves.
\newblock {\em Journal of Computational Technologies}, 10(2):3--14, 2005.

\bibitem{D3S}
Arjen Doelman, Bj{\"o}rn Sandstede, Arnd Scheel, and Guido Schneider.
\newblock The dynamics of modulated wave trains.
\newblock {\em Mem. Amer. Math. Soc.}, 199(934):viii+105, 2009.

\bibitem{Dull_Schneider_NLS}
Wolf-Patrick D{\"u}ll and Guido Schneider.
\newblock Validity of {W}hitham's equations for the modulation of periodic
  traveling waves in the {NLS} equation.
\newblock {\em J. Nonlinear Sci.}, 19(5):453--466, 2009.

\bibitem{Haer}
J{\"o}rg H{\"a}rterich.
\newblock Existence of rollwaves in a viscous shallow water equation.
\newblock In {\em E{QUADIFF} 2003}, pages 511--516. World Sci. Publ.,
  Hackensack, NJ, 2005.

\bibitem{Hwang_Chang}
Shyh~Hong Hwang and Hsueh-Chia Chang.
\newblock Turbulent and inertial roll waves in inclined film flow.
\newblock {\em Phys. Fluids}, 30(5):1259--1268, 1987.

\bibitem{Z_kdv}
Mathew~A. Johnson, Kevin Zumbrun, and Jared~C. Bronski.
\newblock On the modulation equations and stability of periodic generalized
  korteweg-de vries waves via bloch decompositions.
\newblock {\em Physica D: Nonlinear Phenomena}, 239(23-24):2057--2065, 2010.

\bibitem{Jo_Zu_No}
Mathew~A. {Johnson}, Kevin {Zumbrun}, and Pascal {Noble}.
\newblock Nonlinear stability of viscous roll waves.
\newblock {\em ArXiv e-prints, arXiv:1002.0788v1}, 2010.
\newblock Submitted.

\bibitem{Me_Sch}
Ian Melbourne and Guido Schneider.
\newblock Phase dynamics in the complex {G}inzburg-{L}andau equation.
\newblock {\em J. Differential Equations}, 199(1):22--46, 2004.

\bibitem{nd2}
John~H. Merkin and David~J. Needham.
\newblock An infinite period bifurcation arising in roll waves down an open
  inclined channel.
\newblock {\em Proc. Roy. Soc. London Ser. A}, 405(1828):103--116, 1986.

\bibitem{nd}
David~J. Needham and John~H. Merkin.
\newblock On roll waves down an open inclined channel.
\newblock {\em Proceedings of the Royal Society of London. A. Mathematical and
  Physical Sciences}, 394(1807):259--278, 1984.

\bibitem{Noble}
Pascal Noble.
\newblock {\em M\'ethodes de vari\'et\'es invariantes pour les \'equations de
  Saint Venant et les syst\`emes Hamiltoniens discrets}.
\newblock PhD thesis, Universit\'e Paul Sabatier, Toulouse 3, 2003.
\newblock In French.

\bibitem{N2}
Pascal Noble.
\newblock On the spectral stability of roll-waves.
\newblock {\em Indiana Univ. Math. J.}, 55(2):795--848, 2006.

\bibitem{N1}
Pascal Noble.
\newblock Linear stability of viscous roll waves.
\newblock {\em Comm. Partial Differential Equations}, 32(10-12):1681--1713,
  2007.

\bibitem{Schneider_Error}
Guido Schneider.
\newblock Error estimates for the {G}inzburg-{L}andau approximation.
\newblock {\em Z. Angew. Math. Phys.}, 45(3):433--457, 1994.

\bibitem{Serre}
Denis Serre.
\newblock Spectral stability of periodic solutions of viscous conservation
  laws: large wavelength analysis.
\newblock {\em Comm. Partial Differential Equations}, 30(1-3):259--282, 2005.

\bibitem{Whitham}
Gerald~B. Whitham.
\newblock {\em Linear and nonlinear waves}.
\newblock Wiley-Interscience [John Wiley \& Sons], New York, 1974.
\newblock Pure and Applied Mathematics.

\end{thebibliography}

\end{document}